%% file: adjoint_r1.tex
\documentclass[preprint, 10pt, 3p]{elsarticle}

\usepackage{graphicx}
\usepackage{subfig}
\usepackage{amsmath, amssymb}

\journal{Journal of Computational Physics}

\newproof{proof}{Proof}

\newcommand{\hr}{\underline{\rho}}
\newcommand{\hq}{\underline{q}}
\newcommand{\hp}{\underline{\phi}}
\newcommand{\hv}{\underline{v}}
\newcommand{\vf}{\vec{f}}
\newcommand{\hJ}{\underline{\mathcal{J}}}
\newcommand{\pd}[2]{\ensuremath{\frac{\partial{#1}}{\partial{#2}}}}
\newcommand{\dd}[2]{\ensuremath{\frac{d{#1}}{d{#2}}}}

\begin{document}

\begin{frontmatter}

\title{Probability density adjoint for sensitivity analysis of the Mean
of Chaos}

\author[mit]{Patrick J. Blonigan\corref{cor}}
\ead{blonigan@mit.edu}

\author[mit]{Qiqi Wang}
\ead{qiqi@mit.edu}

\address[mit]{Department of Aeronautics and Astronautics,
Massachusetts Institute of Technology, 77 Massachusetts Ave, Cambridge,
MA 02139, United States}


\begin{abstract}
Sensitivity analysis, especially adjoint based sensitivity analysis, is a powerful tool for 
engineering design which allows for the efficient computation of sensitivities with
respect to many parameters. However, these methods break down when used to compute sensitivities
of long-time averaged quantities in chaotic dynamical systems.  

The following paper presents a new method for sensitivity analysis of {\em ergodic} chaotic dynamical systems, the density adjoint method.  The method involves solving the governing equations for the system's invariant measure and its adjoint on the system's attractor manifold rather than in phase-space.  
This new approach is derived for and demonstrated on one-dimensional chaotic maps and the three-dimensional Lorenz system.  It is found that the density adjoint computes very finely detailed adjoint distributions and accurate sensitivities, but suffers from large computational costs.  

\end{abstract}

\end{frontmatter}

\section{Introduction}
\label{s:introduction}

Sensitivity analysis of systems governed by ordinary differential
equations and partial differential equations are important in many
fields of science and engineering.  Its goal is to compute sensitivity
derivatives of key quantities of interest to parameters that influence
the system.  Applications of sensitivity analysis in science and
engineering include design optimization, inverse problems,
data assimilation, and uncertainty quantification.

Adjoint based sensitivity analysis is especially powerful in many
applications, due to its efficiency when the number of parameters is
large.  In airplane design, for example, the number of geometric parameters that
define the aerodynamic shape is very large.  As a result, the adjoint
method of sensitivity analysis proved very successful for aircraft design
 \cite{Jameson:1988:adj}.  
Similarly, the
adjoint method has been an essential tool for solving inverse problems in
seismology, and for assimilating observation data for weather
forecasting.

Sensitivity analysis for chaotic dynamical systems is important because of the 
prevalence of chaos in many scientific and engineering fields.  One example is 
chaotic aero-elastic oscillations of aircraft wings and control surfaces.  
In this example, and in other applications with periodic or 
chaotic characteristics, statistical averaged quantities such as mean 
stresses and mean aerodynamic forces are of interest.  Therefore, the 
general problem this paper seeks a solution to is:
 
\begin{equation}
\mbox{Given } \frac{d\vec{x}}{dt} = \vec{f}(\vec{x}, \xi), \quad
\overline{J} = \lim_{T\rightarrow\infty} \frac1T \int_0^T J(\vec{x},\xi)
dt,\quad
\mbox{Compute } \frac{\partial \overline{J}}{\partial \xi}
\label{e:gen_problem}
\end{equation}

Sensitivity analysis for chaotic dynamical systems is difficult
because of their sensitivity to the initial condition, known as the "Butterfly 
Effect".  Slightly different initial conditions will result in very different 
solutions, which diverge exponentially with time \cite{Lorenz:1963:det}.  
This also results in exponential growth of sensitivities and therefore the 
sensitivity of long-time averaged quantities is not equal to the long-time 
average sensitivities of chaotic systems \cite{Lea:2000:climate_sens}.  Because 
the derivative and long-time average do not commute, the traditional adjoint 
method computes sensitivities that diverge, as shown in the work done by 
Lea et al. \cite{Lea:2000:climate_sens}.  

Prior work in this area includes the ensemble-adjoint method proposed by Lea et 
al. \cite{Lea:2000:climate_sens} and then applied to an ocean circulation model 
with some success \cite{Lea:2002:ocean}.  Eyink et al. went on to generalize 
the method \cite{Eyink:2004:ensmbl}.  The ensemble-adjoint involves averaging 
over a large number of ensemble calculations and the resulting high 
computational costs make this method intractable for many applications. 

Climate sensitivity analysis of chaotic systems based on the probability
density function in phase space has a long history.
A perturbation to the dynamical system causes a corresponding perturbation in
the stationary probability density function.  This correspondence is
governed by the Fokker-Planck equation, also known as the Liouville equation
for conservative dynamical systems.  An analysis based on the Fokker-Planck
equation produces the Fluctuation Dissipation Theorem
\cite{PhysRev.32.110, 0034-4885-29-1-306}.  For conservative and nearly
conservative dynamical systems, the Fluctuation Dissipation Theorem
can be used to accurately compute climate sensitivities \cite{leith1975climate}.
Several improved algorithms based on Fluctuation Dissipation Theorem have since
been developed for computing climate sensitivity of non-conservative systems
\cite{majda2005information, 0951-7715-20-12-004}.  In particular, an approach
based on numerically solving the Fokker-Planck equation has been demonstrated for
strongly dissipative chaotic dynamical systems \cite{Thuburn:2005:FP}.
This approach involves 
finding a probability density function which satisfies a Fokker-Plank equation 
to model the climate.  The adjoint of this Fokker-Planck equation is then used 
to compute derivatives with respect to statistically averaged quantities.  
However, their method requires adding diffusion into the system,
potentially making the computed sensitivity inaccurate.

This paper presents a new method for computing sensitivity of
mean quantities in ergodic chaotic dynamical systems based on a Fokker-Planck type 
formulation.  The key idea is to describe the objective function $\overline{J}$ as 
an average in phase space as in \cite{Thuburn:2005:FP}:
\begin{equation}
\overline{J} = \int_{R^n} J(\vec{x}) \rho_s(\vec{x}) d\vec{x}
\label{e:phase_space_avg}
\end{equation}
The density of the invariant measure of the chaotic system $\rho_s(\vec{x})$ is 
governed
by a probability density equation, whose adjoint equation can be solved to
compute the desired sensitivities.  


Our method, the density adjoint method, relies on the following assumptions:

\begin{itemize}
\item The chaotic dynamical system is ergodic
\item The system has a smooth stationary density distribution $\rho_s$.      
\item Perturbations to long time averaged quantities of interest $\overline{J}$ depend mainly on perturbations to $\rho_s$ on the attractor surface and less so on perturbations to the position and shape of the attractor manifold.  
\end{itemize}

The rest of this paper is organized as follows: Section \ref{s:smoothness}
discusses the well-posedness of the problem (equation \eqref{e:gen_problem}) by analyzing 
the
differentiability of the time averaged quantities $\overline{J}$ for a few discrete and continuous chaotic dynamical systems.
Section \ref{s:density1d} presents the probability density adjoint method for chaotic,
1D iterated maps.  Section \ref{s:densitylorenz} extends our method to continuous
dynamical systems, with the Lorenz attractor as an example.  Section \ref{s:densitylorenz} also includes some considerations for minimizing errors in the density adjoint method and discusses the limitations of the method.  Section \ref{s:conclusion} concludes this paper.

\section{Smoothness of the Mean and Stationary Density
Distribution of Chaos}
\label{s:smoothness}

Not every chaotic dynamical systems has differentiable mean quantities
$\overline{J}$.   Hyperbolic systems, a class of dynamical systems with
ideal attractors, are known to have mean quantities that respond
differentiably to small perturbations in its parameters \cite{Bonatti:2005:Hyper}.
Chaotic systems whose mean quantities are
differentiable to perturbations are generally classified as
quasi-hyperbolic systems \cite{Bonatti:2005:Hyper}.  Other chaotic dynamical systems 
are known as non-hyperbolic.  In these non-hyperbolic systems, the mean quantities
are usually not differentiable, or even continuous as the parameters
vary.  In fact,
the long time average for non-hyperbolic systems may have nontrivial
dependence on the initial condition, indicating that the mean quantity is not
even well-defined.  

The class of the system can be related to properties of the stationary density distribution of a system, $\rho_s(\vec{x})$.  The stationary density distribution is a density distribution in phase space that is invariant under the dynamical system.  For
hyperbolic and quasi-hyperbolic systems, it can be rigorously
characterized as the Sinai-Ruelle-Bowen (SRB) measure \cite{Rulle:1997:SRB}.
The stationary density can be computed by evolving the
dynamical system, with the initial condition drawn from an arbitrary, continuous density
distribution in phase space.  

It can be shown that systems with smooth density distributions have differentiable mean quantities\footnote{However, hyperbolic/quasi-hyperbolic systems do not necessarily have smooth density distributions}.  Denote the governing equations for the stationary density distribution, $\rho_s$ as:

\begin{equation}
 \mathcal{L} \rho_s = 0
 \label{e:density_gov}
\end{equation}

It will be shown that $\mathcal{L}$ is related to the Frobenius-Perron operator\footnote{$\mathcal{L} \rho_s = P \rho_s - \rho_s$, where $P$ is the Frobenius-Perron operator.} for 1D maps in section \ref{s:density1d} or the Fokker-Planck equation in section \ref{s:densitylorenz}.  The sensitivity of a mean quantity, $\overline{J}$, to some parameter $\xi$ can be expressed as a function of the stationary density sensitivity, $\pd{\rho_s}{\xi}$, by differentiating equation \eqref{e:phase_space_avg}:    

\begin{equation}
\dd{\overline{J}}{\xi} = \int_{R^n} J(\vec{x}) \pd{\rho_s}{\xi} \ d\vec{x}
\label{e:phase_space_sens}
\end{equation}

\noindent where $\pd{\rho_s}{\xi}$ satisfies the linearization of equation \eqref{e:density_gov}:

\begin{equation}
 \mathcal{L} \pd{\rho_s}{\xi} = -\left(\pd{ }{\xi}\mathcal{L}\right) \rho_s
 \label{e:lin_density}
\end{equation}

The operator on the right hand side, $\pd{ }{\xi}\mathcal{L}$, can be shown to include spatial derivatives (see section \ref{ss:grad1d} and \ref{ap:NDcont} for the Frobenius-Perron operator and the Fokker-Planck equation, respectively).  Therefore, if $\rho_s$ is differentiable in phase space and $\mathcal{L}$ is not poorly conditioned, $\frac{\partial \rho_s}{\partial \xi}$ is finite according to equation \eqref{e:lin_density}.  If this is the case, then from equation \eqref{e:phase_space_sens}, \dd{\overline{J}}{\xi} is also finite, so $\overline{J}$ is differentiable.  

\begin{figure}[htbp]
  \centering
  \begin{minipage}[b]{3.125in}
    \includegraphics[width=3.2in]{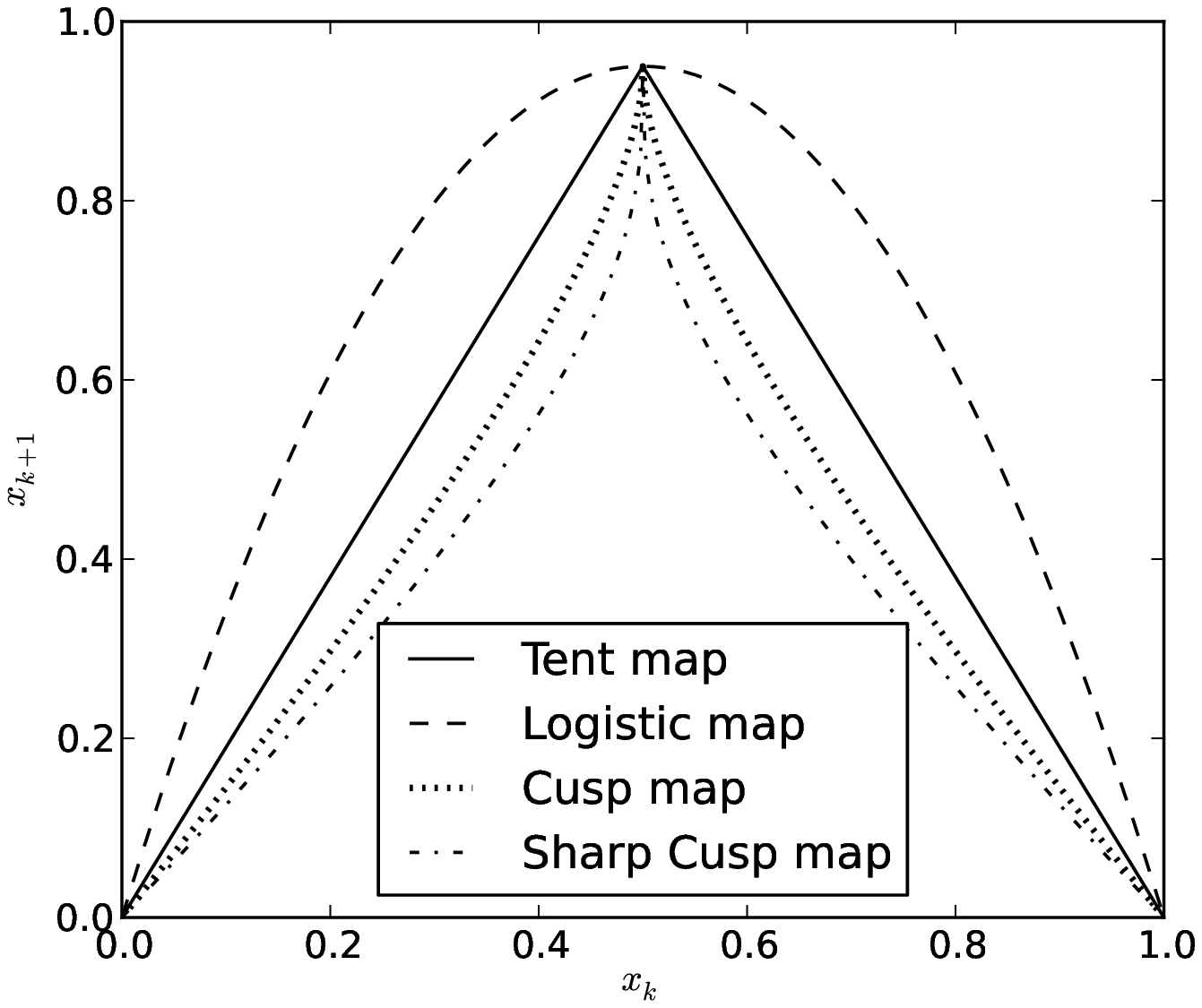}
    \caption{The shape of the logistic, tent and Cusp maps.}
    \label{f:maps}
  \end{minipage}
  \hspace{0.125 in}
  \begin{minipage}[b]{3.125in}
    \includegraphics[width=3.2in]{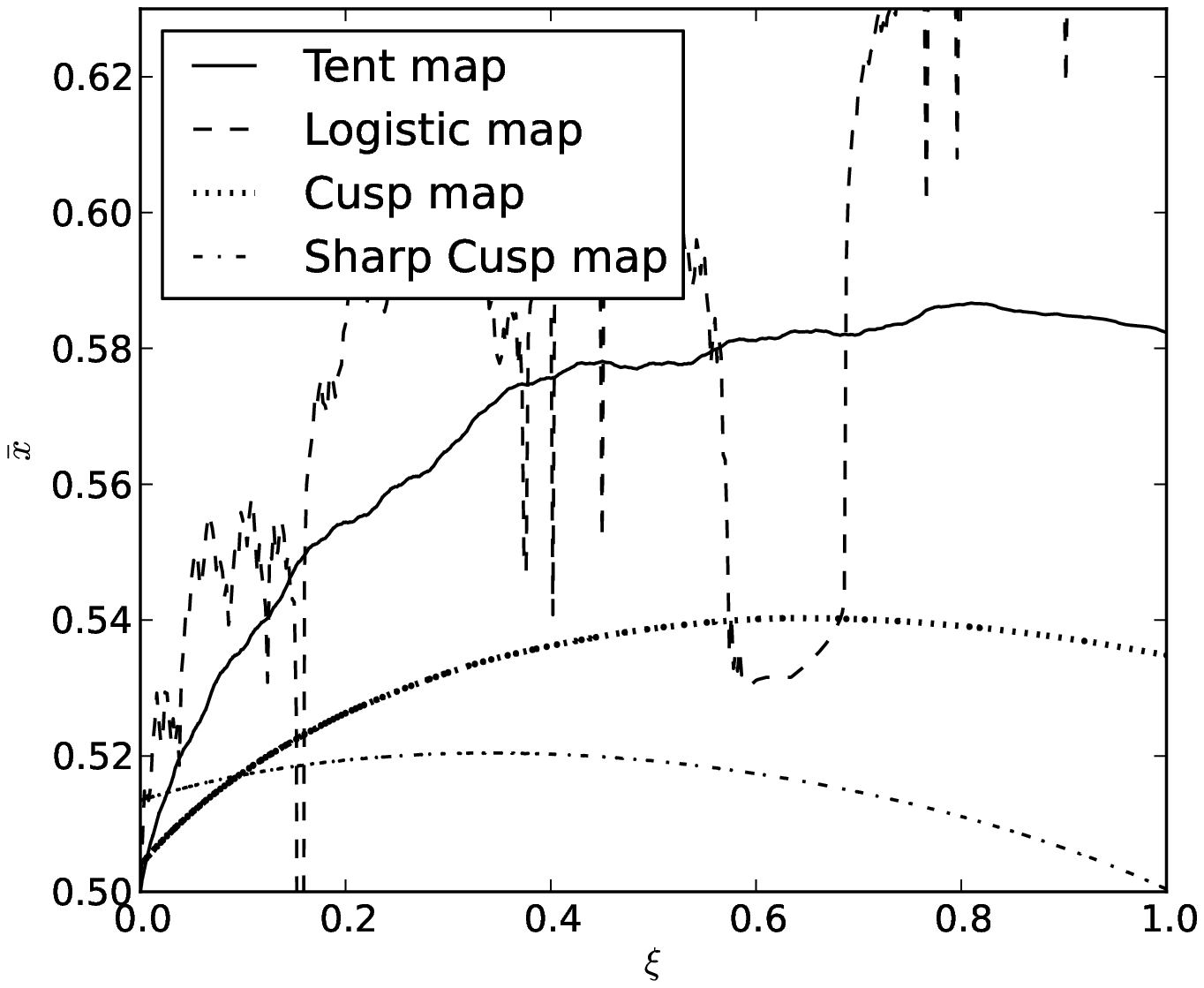}
    \caption{Smoothness of $\overline{x}$ as a function of parameter $\xi$.  }
    \label{f:smoothness}
  \end{minipage}
\end{figure}

The relation between the smoothness of $\rho_s$ and the differentiability of $\overline{J}$ can be demonstrated numerically as well.  We first study three parametrized 1D chaotic maps:
\begin{enumerate}
\item The logistic map
\begin{equation}
x_{k+1} = F_{logistic}(x_k) = \left(4-\frac{\xi}{4}\right) x_k (1-x_k)
\end{equation}
\item The tent map
\begin{equation} \label{tentmap}
x_{k+1} = F_{tent}(x_k) = \left(2-\frac{\xi}{2}\right) \min(x_k, 1-x_k)
\end{equation}
\item The Cusp map
\begin{equation} \label{cuspmap}
x_{k+1} = F_{cusp}(x_k) = \left(1-\frac{\xi}{4}\right)
          \left(1 - \left|\frac12-x_k\right|
                  - \sqrt{\left|\frac14-\frac{x_k}2\right|}\right)
\end{equation}
\noindent We also consider another ``sharp'' version of the Cusp map
\begin{equation}
x_{k+1} = F_{cusp}(x_k) = \left(1-\frac{\xi}{4}\right)
          \left(1 - \left|\frac12-x_k\right|
                  - \left(\left|\frac14-\frac{x_k}2\right|\right)^{0.3}\right)
\end{equation}
\end{enumerate}
In all 3 maps, the parameter $\xi$ controls the height of the maps.
Figure \ref{f:maps} shows the logistic map for $\xi=0.8$, the tent map for
$\xi=0.2$ and the two Cusp maps for $\xi=0.2$.

Although the logistic map, the tent map and the Cusp map have the same monotonic 
trends in $[0, 0.5]$ and $[0.5,1]$, the smoothness of their mean
\begin{equation}
\overline{x} = \lim_{N\rightarrow\infty}\frac1N \sum_{k=1}^N x_k
\end{equation}
with respect to the parameter $\xi$ are very different.
Figure \ref{f:smoothness} plots the mean $\overline{x}$ of the three
chaotic maps against the parameter $\xi$.  The mean is approximated as
\begin{equation}
\overline{x} \approx \frac1{N\,M} \sum_{i=1}^{M} \sum_{k=n_0}^{N+n_0} x_{i,k},
\quad x_{i,k+1} = F(x_{i,k}),\quad
\quad M=1000, N=50000, n_0=1000,
\end{equation}
and $x_{i,0}$ are uniformly randomly sampled in $[0.25,0.75)$ so that $\overline{x}$ is computed from $M$ different trajectories.  $n_0$ is the number of ``spin up'' iterations for the map, to eliminate any transient features of the trajectories.  

The mean $\overline{x}$ of the logistic map appears to be discontinuous
with respect to $\xi$.
The mean of the tent map appears to be continuous and differentiable
with respect to $\xi$, but it is difficult to assess its higher order
smoothness.  The mean of the two Cusp maps appears to be smoother than
the tent map.

\begin{figure}[htb!] \centering
\includegraphics[width=3.2in]{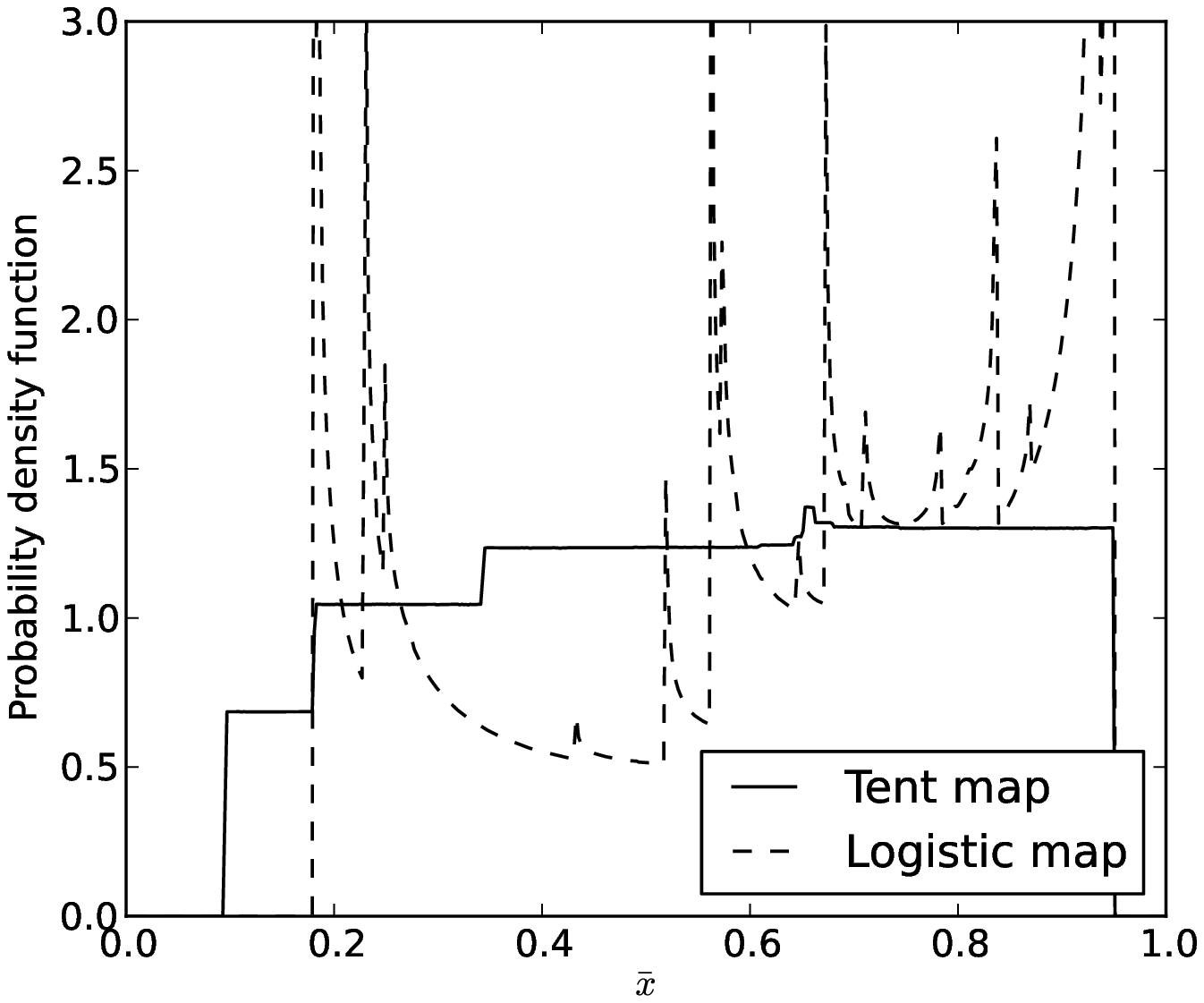}
\includegraphics[width=3.2in]{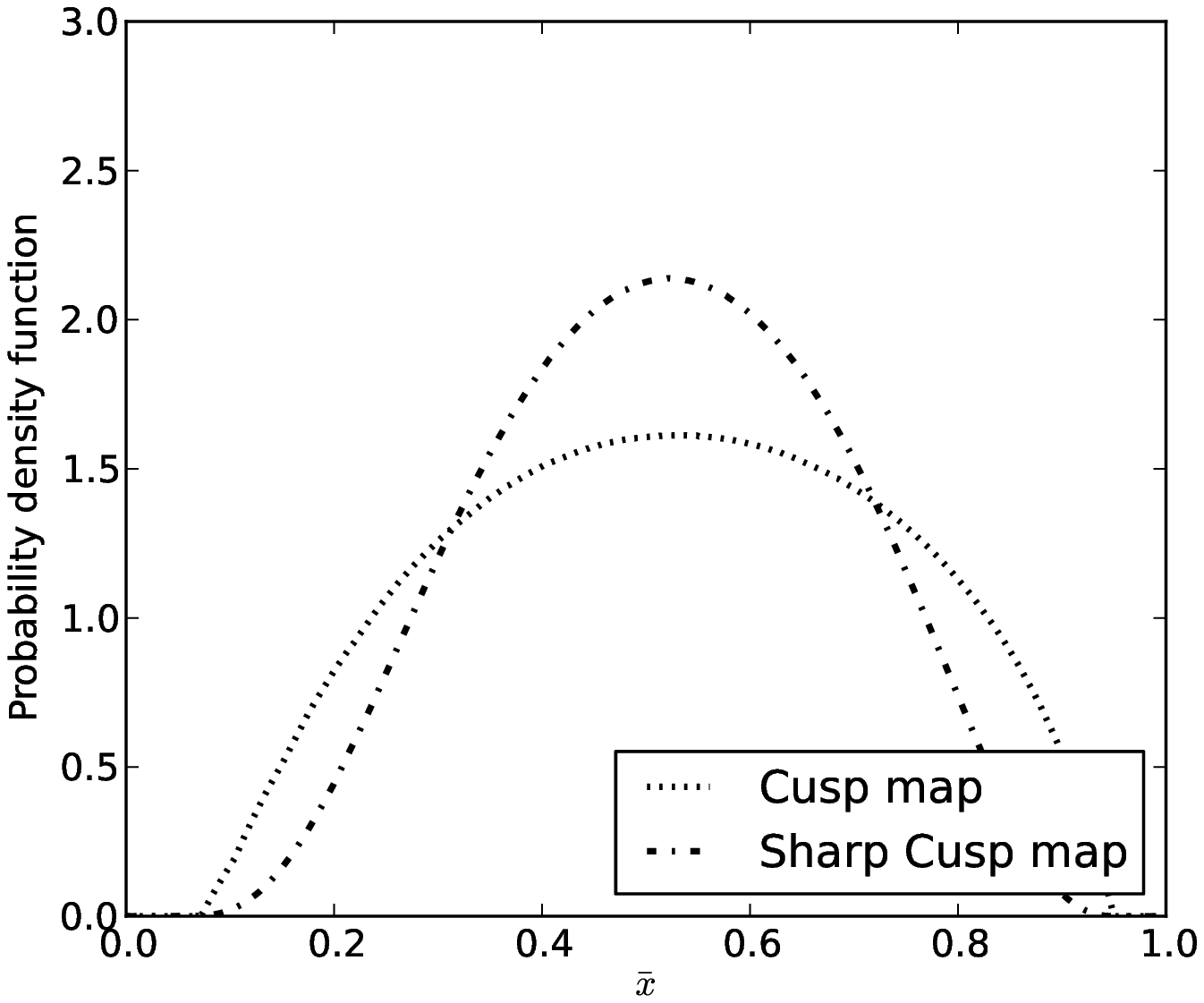}
\caption{Stationary density}
\label{f:density}
\end{figure}

Figure \ref{f:density} shows the stationary density distribution of the
logistic map, the tent map and the two Cusp maps.  The logistic map
has a stationary density function that concentrates at discrete
points, as is evident from the peaks in the density function.  The
stationary density function of the tent map is bounded, but appears to
contain discontinuities.  The density of the Cusp map is continuous;
while the density of the sharp Cusp map appears to be the most smooth.
We find that the maps with smoother mean quantities tend to have smoother 
stationary density distributions.  

The same conclusion can be drawn for continuous dynamical systems.
Here, we analyze the mean quantities and the stationary density 
distributions of the two most well known chaotic attractors:
the R\"{o}ssler attractor
\begin{equation}
\frac{dx}{dt} = -y - z, \quad
\frac{dy}{dt} = x + a\,y, \quad
\frac{dz}{dt} = b + z(x-c)
\end{equation}
and the Lorenz attractor
\begin{equation}
\frac{dx}{dt} = s (y - x), \quad
\frac{dy}{dt} = x (r - z) - y, \quad
\frac{dz}{dt} = x\,y - b\,z\;.
\label{e:lorenz}
\end{equation}
For the R\"{o}ssler attractor, we analyze how
$\overline{x}$ and $\overline{z}$ change as the
parameter $c$ varies.  For the Lorenz attractor, we know that
$\overline{x} \equiv \overline{y} \equiv 0$
due to symmetry of the governing equation.  Therefore,
we focus on the nontrivial quantities
$\overline{z}$ and $\overline{x^2}$
as the Rayleigh number $r$ varies.

\begin{figure}[htbp]
  \centering
  \begin{minipage}[b]{3.125in}
    \includegraphics[width=3.2in]{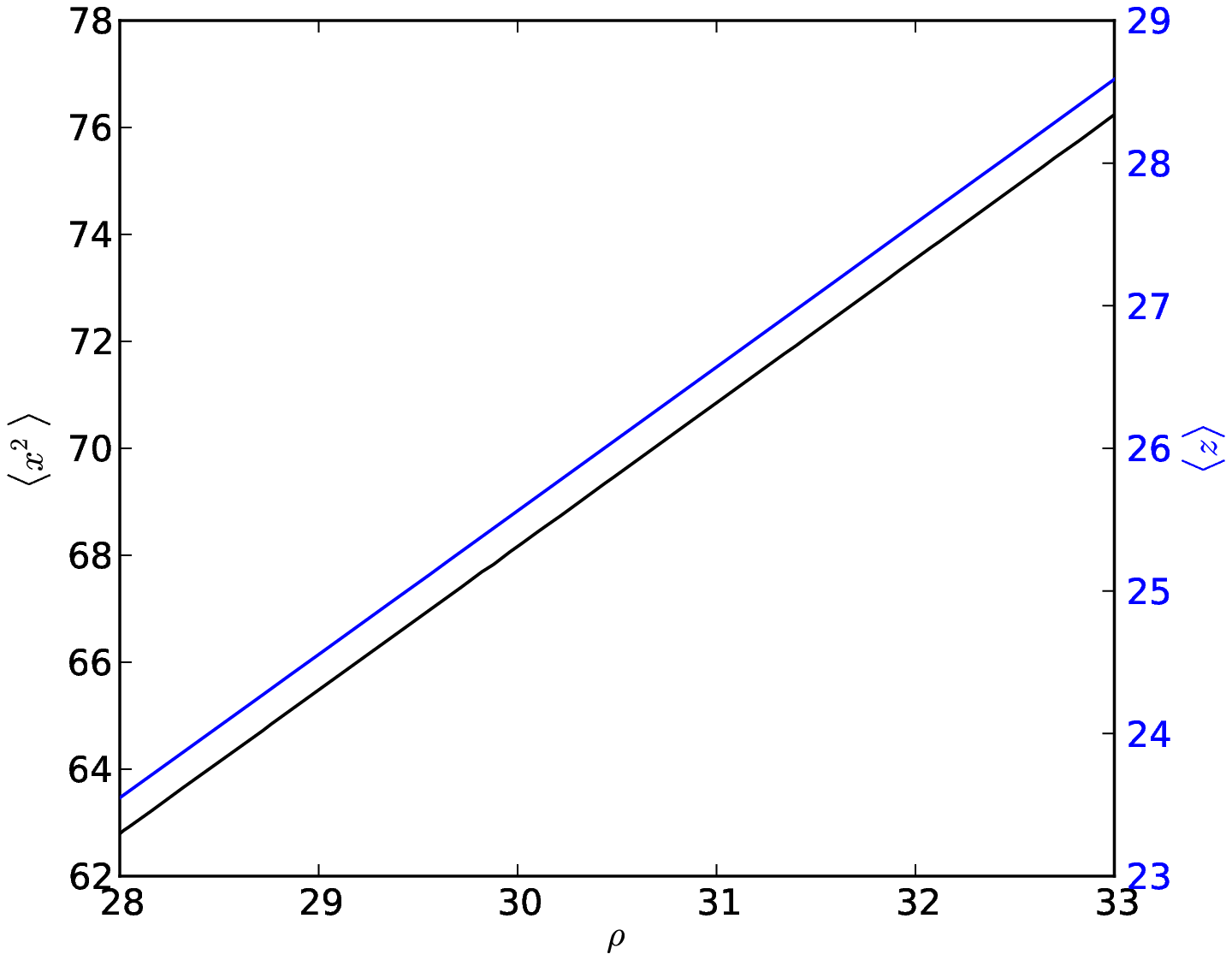}
    \caption{$\overline{z}$ and $\overline{x^2}$ of the Lorenz attractor
             as $r$ varies.}
    \label{f:lorenzPoincareMap}
  \end{minipage}
  \hspace{0.125 in}
  \begin{minipage}[b]{3.125in}
    \includegraphics[width=3.2in]{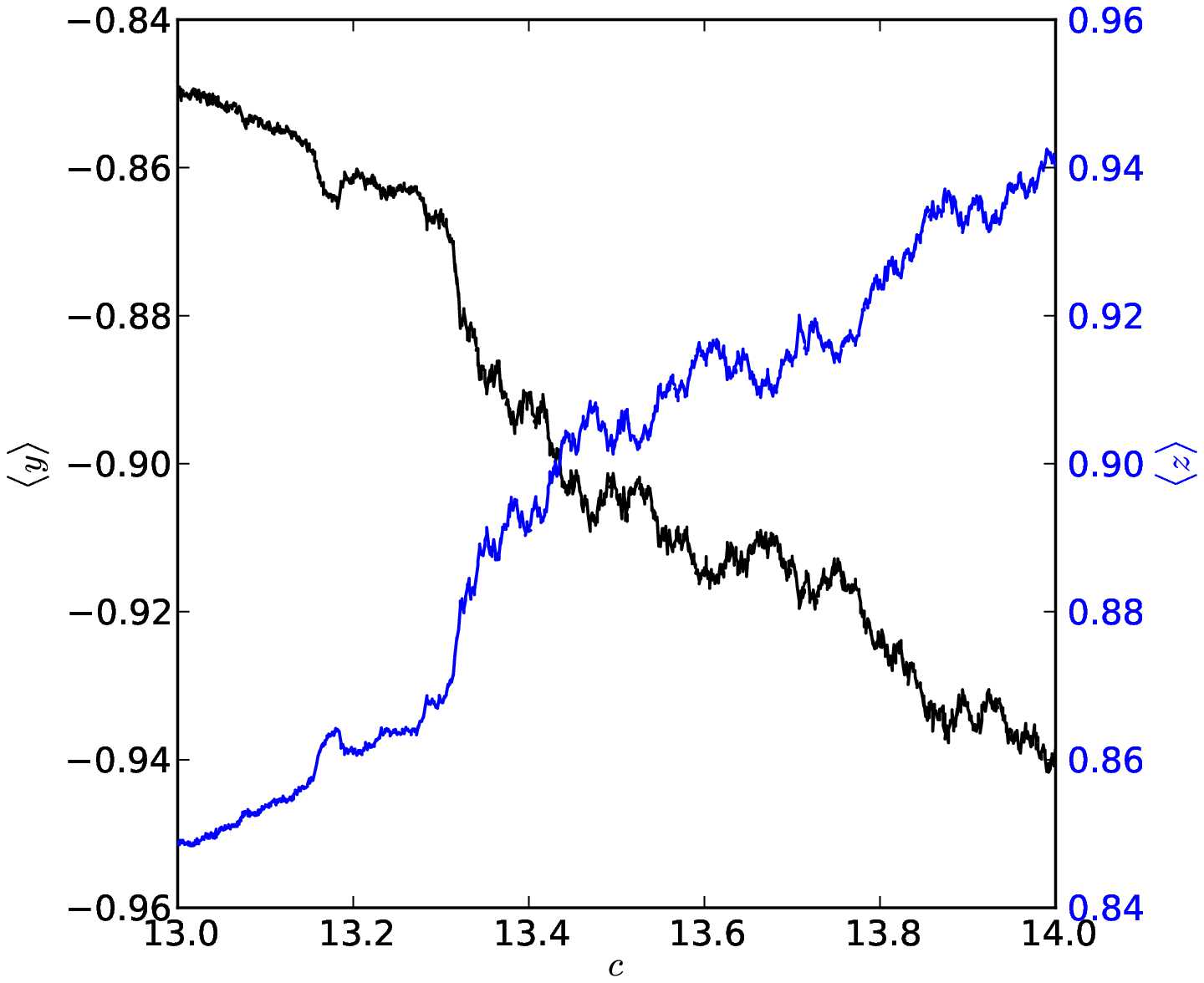}
    \caption{$\overline{x}$ and $\overline{z}$ of the R\"{o}ssler attractor
             as $c$ varies.}
    \label{f:rosslerPoincareMap}
  \end{minipage}
\end{figure}

Figures \ref{f:lorenzPoincareMap} and \ref{f:rosslerPoincareMap} show
how the mean quantities respond to parameter changes for the Lorenz
attractor and the R\"{o}ssler attractor.  The R\"{o}ssler attractor has similar
behavior to the logistic map.  The mean quantities are not smooth
functions of its parameter $c$.  The Lorenz attractor has mean quantities
that are smooth functions of its parameter $r$.

\begin{figure}[htb!] \centering
\subfloat[Lorenz attractor at $r=28,
s=10,b=8/3$.]
{\includegraphics[width=5.2in,trim=0 10 0 25,clip]{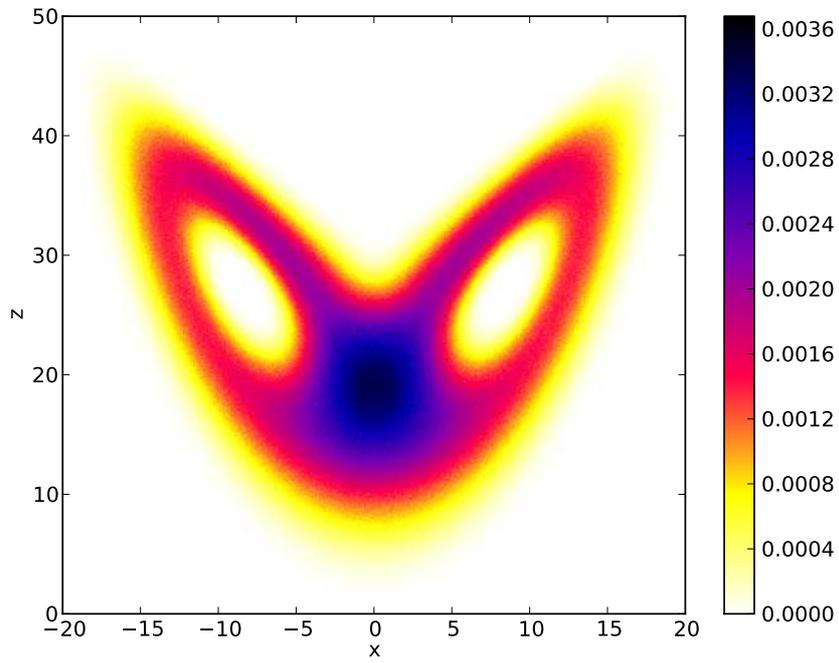}}\\
\subfloat[R\"{o}ssler attractor at $a=b=0.1, c=14$.]
{\includegraphics[width=5.2in,trim=0 10 0 25,clip]{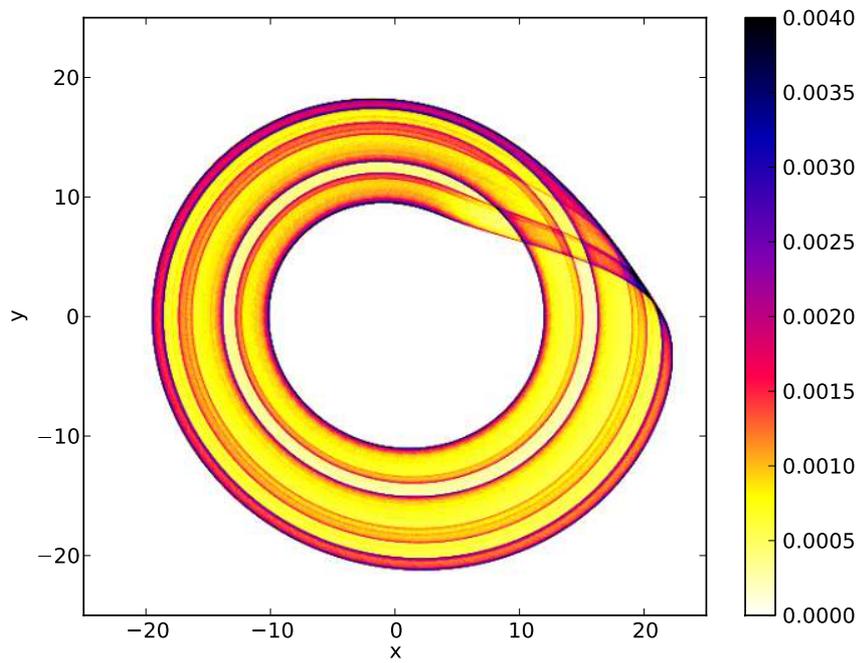}}
\caption{The stationary density of the Lorenz attractor and the R\"{o}ssler
attractor projected into the $xz$ and $xy$ planes, respectively.}
\label{f:lorenzrosslerdensity}
\end{figure}

Figures \ref{f:lorenzrosslerdensity} a and b show the stationary density distributions 
projected onto the $xz$ and $xy$ planes respectively for the Lorenz and R\"{o}ssler
attractors.  As was the case for the 1D maps, the density of the Lorenz attractor, whose
mean quantities vary smoothly with respect to parameter changes has a smooth density 
distribution.  A number of discontinuities are present in the density distribution 
of the R\"{o}ssler attractor, whose mean quantities do not exibit smooth variation with 
respect to parameter changes.  

The relationship between the smoothness of stationary density distribution and the smoothness of mean quantities provides a justification for the method
developed in this paper.  If the stationary
density distribution in phase space is smooth on its attractor 
manifold, the mean quantities are differentiable with
respect to the parameters of a chaotic dynamical system.  
The density distribution can be accurately solved by 
discretizing its governing equation, the Fokker-Planck equation, 
on its attractor manifold. Sensitivity derivatives of the mean quantities 
with respect to system parameters can then be computed via sensitivity 
analysis of the Fokker-Planck equation.

\section{Density adjoint for chaotic 1D maps}
\label{s:density1d}

This section uses the parameterized cusp map as an example
to illustrate the density adjoint method.  This 1D map is
defined as

\begin{equation}
x_{k+1} = F_{cusp}(x_k) = 1 - \xi|2x -1| - (1 - \xi) \sqrt{|2x - 1|}
\label{e:cusp}
\end{equation}

Where the parameter $0\le\xi\le 1$ defines the shape of the map.
When $\xi = 1$, the map is a tent map (\ref{tentmap}); when $\xi=0.5$, the
map is a cusp map (equation (\ref{cuspmap})).
The density adjoint method was used to compute the sensitivity of the 
mean $\bar{x}$ with respect to the parameter $\xi$.

\subsection{Computing Stationary Density}
\label{ss:stat_dist}

The stationary density distribution $\rho_s(x)$
is a one dimensional probability density distribution determined by a given mapping 
function $x_{k+1} = F(x_k)$.  It is governed by the Frobenius-Perron equation \cite{Ding:1991:Frob-Per}, and defines the 
probability that an initial point $x_0$ will be mapped to some region $\Delta 
x$ after infinitely many mappings. Consider a series of random variables $X_k$ satisfying
$X_{k+1} = F(X_k)$ for all $k\ge0$.  The distribution of $X_k$ converges to the stationary 
distribution as $k\to \infty$ whenever $X_0$ has a finite distribution function.  
Denote the Frobenius-Perron operator $P$ 
as the map from the probability distribution $\rho_k$ of $X_k$ to the probability distribution 
$\rho_{k+1}$ of $X_{k+1}$ \cite{Ding:1991:Frob-Per}.  Then $\rho_s = \lim_{k \to \infty} P^k(\rho_0)$ for any finite $\rho_0$.  
An equivalent statement is that $\rho_s(x)$ is an eigenfunction of the operator $P$, with 
an eigenvalue of one:

\begin{equation}
(P\rho_s)(x) = \rho_s(x) \ , \quad x\in[0,1]
\label{e:eigen_problem}
\end{equation}

\noindent The operator $P$ in equation (\ref{e:eigen_problem}) is the Frobenius-Perron operator 
defined in \cite{Ding:1991:Frob-Per} as: 

\begin{equation}
\int_0^1 P \rho(x) dx = \int_0^1 \rho(F(x)) dx
\label{e:FP}
\end{equation}

\noindent To derive $P$, recall that probability density is conserved in our domain, 
phase space, by the normalization axiom of probability.  

\begin{figure}[htb!] \centering
\includegraphics[width=3.2in]{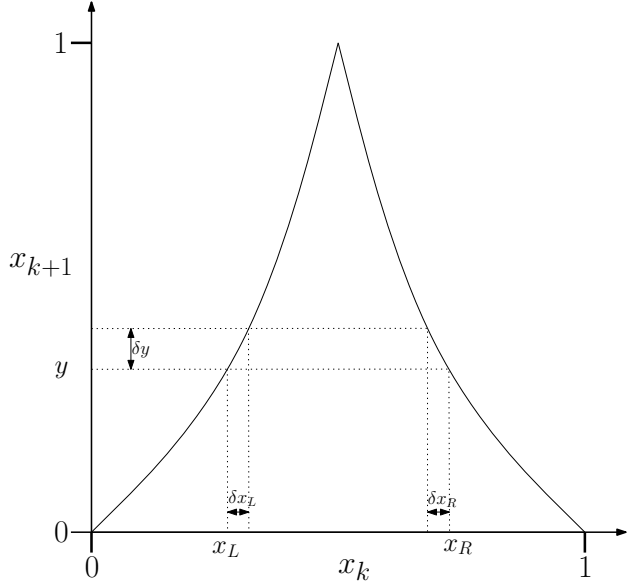}
\caption{Density Mapping for the cusp map}
\label{f:density_cons}
\end{figure}

\noindent In the case of the map shown in figure \ref{f:density_cons}, the integral of 
the density contained in the small intervals $\delta x_L$ and $\delta 
x_R$ will be mapped into the interval $\delta y$.  This can be written as 
follows, where $y=F(x_L)=F(x_R)$: 

\begin{equation*}
\int_y^{y + \delta y} \rho_{k+1}(s) \ ds = \int_{x_L}^{x_L + \delta x_L} \rho_k(s) \ ds + 
\int_{x_R}^{x_R + \delta x_R} \rho_k(s)\ ds
\end{equation*}

\noindent Differentiating with respect to $s$ and dividing both sides by $dy/ds$, 
an expression for the mapping of density is obtained: 

\begin{equation}
\rho_{k+1}(y) = \frac{1}{|F'(x_L)|}\rho_k(x_L) + 
\frac{1}{|F'(x_R)|}\rho_k(x_R)
\label{e:density_mapping}
\end{equation}

\noindent Where $F'(x) = \frac{dF}{dx}$. Ding and Li \cite{Ding:1991:Frob-Per} compute $\rho_s$ by using finite elements 
to construct a discrete approximate of $P$.  Both linear and higher order elements 
were investigated and $\rho_s$ was correctly computed for a number of 1D maps including 
the tent map.  

We construct a finite difference discretization of the Frobenius-Perron operator $P$ 
based on equation (\ref{e:density_mapping}).  The interval $[0,1]$
is discretized into $n$ equally spaced nodes, with $y_i = \frac{i-1}{n-1}$.  
We represent the discretized 
version of the linear operator $P$ as an $n$ by $n$ matrix $P_n$. From equation
 (\ref{e:density_mapping}), for the discretized density distributions 
$\hr_k \equiv (\rho_k(y_1),\rho_k(y_2),...,\rho_k(y_n))$
and $\hr_{k+1} \equiv (\rho_{k+1}(y_1),\rho_{k+1}(y_2),...,\rho_{k+1}(y_n))$:

\[
 \hr_{k+1} =  P_n \hr_k
\]

\begin{figure}[htbp]
  \centering
  \begin{minipage}[b]{3.125in}
    \includegraphics[width=3.2in]{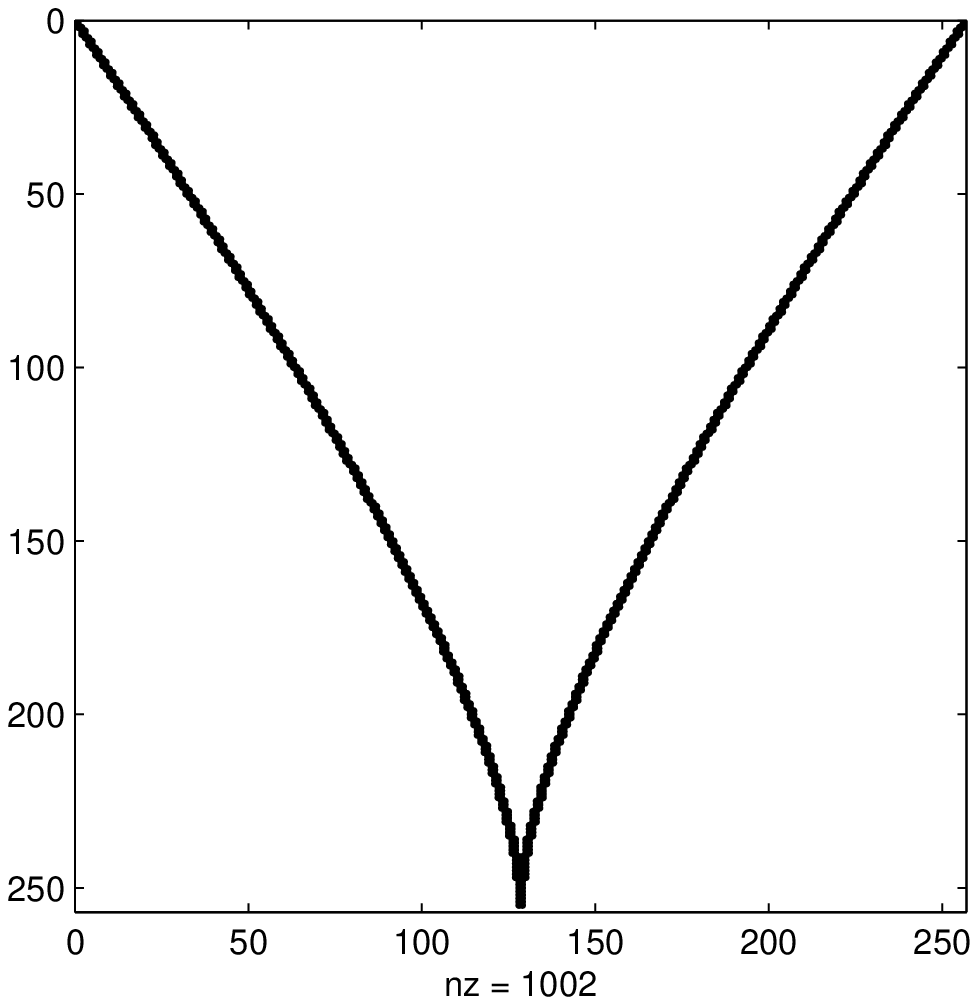}
    \caption{Cusp map transition matrix $P_n$ structure for $\xi=0.5$. }
    \label{f:cusp_matrix}
  \end{minipage}
  \hspace{0.125 in}
  \begin{minipage}[b]{3.125in}
    \includegraphics[width=3.2in]{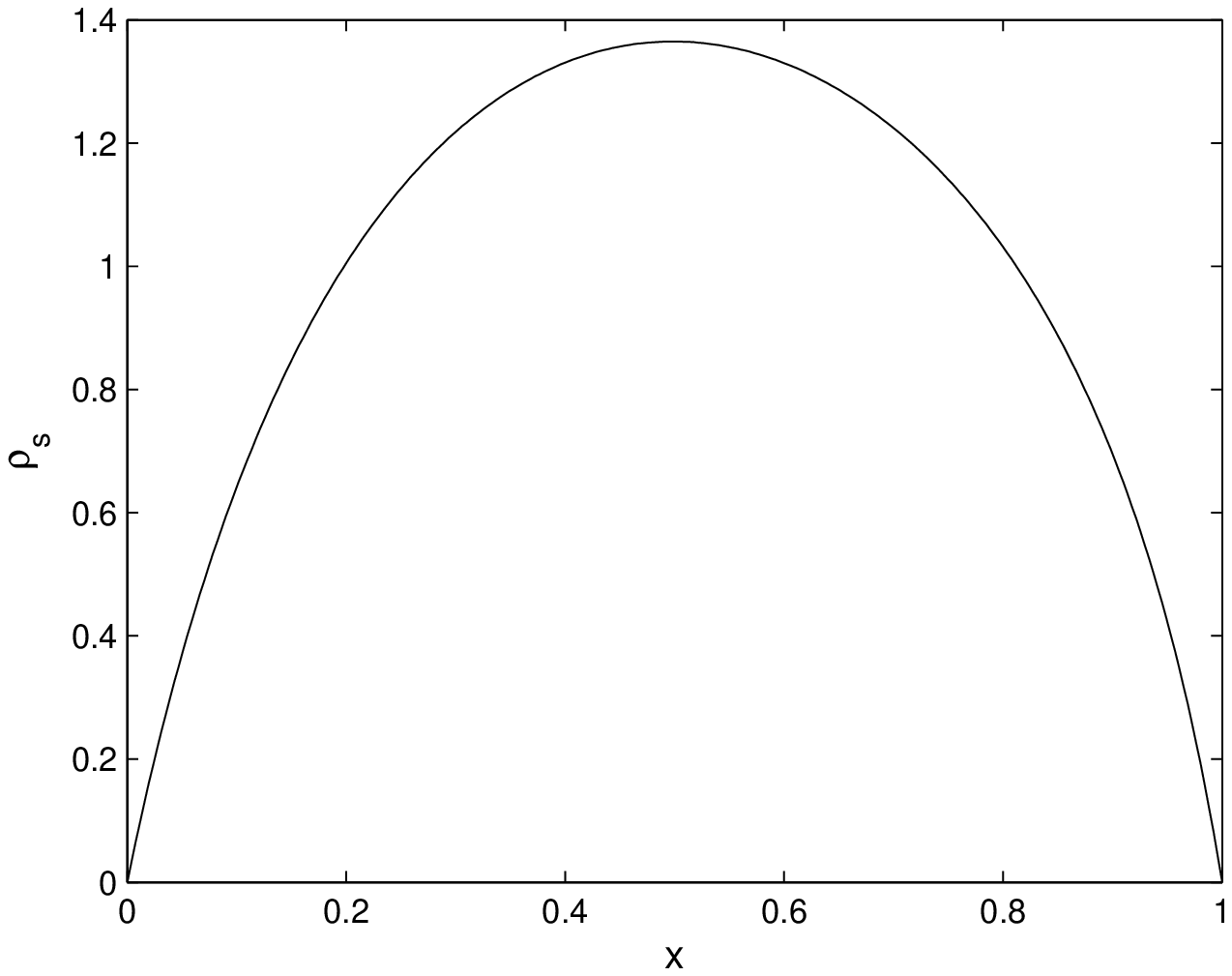}
    \caption{Cusp map density distribution $\rho_s$ for $\xi=0.5$.  Generated 
with 256 nodes. }
    \label{f:cusp_density}
  \end{minipage}
\end{figure}

The matrix $P_n$ is constructed by finding $x_{Li}, x_{Ri} = F^{-1}(y_i)$.  This is done by 
computing the inverse functions associated the left and right sides of $F(x)$ with Newton's 
method. Next, $F'(x)$ is determined at all $x_{Li}$ and $x_{Ri}$. 
In most cases, $x_{Li}, x_{Ri}$ will not be equal to any $y_k$ from the discretization.  
To account for this, $\rho(x_{Li})$ and $\rho(x_{Ri})$ 
are found by linear interpolation between the two nearest nodes.  This means that each row 
of $P_n$ will typically contain two pairs of non-zero entries, one pair for the 
right side of $F(x)$, the other for the left side.  For a uniform discretization of $y_i$, the 
non-zero entries will form the shape of $F(x)$ upside down in the matrix, as 
shown in figure \ref{f:cusp_matrix}.  It is important to note that although 
(\ref{e:density_mapping}) is derived assuming conservation of probability mass, 
$P_n$ does not conserve probability mass.  Unlike $P$, the largest eigenvalue 
$\lambda$ of $P_n$ is not exactly one, due to numerical error from the interpolation.  
To use $P_n$ to compute $\rho_s(x)$ with a power iteration, $\rho_s(x)$ must be 
scaled after each iteration such that its integral is equal to one.  

As shown by figure \ref{f:cusp_density}, $\rho_s(x)$ for the cusp map is continuous, 
showing that the objective function is continuous with respect to $\xi$ and the 
sensitivity with respect to $\xi$ is defined.

\subsection{Computing gradients using the density adjoint}
\label{ss:grad1d}

By the definition of the Frobenius-Perron operator $P$, a perturbation to the mapping function $F(x)$ leads to a 
perturbation to $P$. As $\rho_s(x)$ is the first eigenfunction of the Frobenius-Perron 
operator, there is a density perturbation $\delta\rho_s(x)$ corresponding to a perturbation
to the operator. A perturbation to a mean quantity $\delta\bar{J}$ can be computed 
from $\delta\rho_s(x)$ using the following expression, where $J(x)$ is the quantity of 
interest:

\begin{equation}
\delta\bar{J} = \int_0^1 J(x) \ \delta \rho_s(x) \ dx
\label{e:FP_fwd}
\end{equation}

\noindent $\delta\bar{J}$ can also be computed using the adjoint density $\phi$

\begin{equation}
 \delta \overline{J} = \int_0^1 \phi(x) \ \delta P \ \rho_s(x) \ dx
 \label{e:FP_bkw}
\end{equation}

\noindent where $\phi$ satisfies the adjoint equation:

\begin{equation}
 P^*\phi - \lambda \phi = \overline{J} - J
 \label{e:adj_cont}
\end{equation}

\noindent {\em For a more detailed derivation of the adjoint equation, see appendix \ref{ap:1Dcont}} \\

\noindent $\lambda$ is the first eigenvalue of the operator and is equal to one.  

\noindent The term $\delta P \rho_s$ in equation 
\ref{e:FP_bkw} can be found by considering the mapping of probability mass.  From equation (\ref{e:eigen_problem}):

\[
 \rho_s + \delta P \rho_s = (P + \delta P) \rho_s 
\]

\noindent Assuming a small perturbation $\delta P$ (and therefore a small $\delta F$):  

\begin{equation}
 \int_0^y \delta P \rho_s \ ds = \rho_s(x_L) \delta x_L - \rho_s(x_R) \delta x_R 
 \label{e:cdf_perturb2}
\end{equation}

\noindent For a small perturbation $\delta F$, it can be shown that (see appendix \ref{ap:1Dgrad}):  

\[
 \frac{\delta F}{\delta x} \approx F'(F^{-1}(y))
\]

\noindent Substituting into equation (\ref{e:cdf_perturb2}) and differentiating with respect 
to $y$:

\begin{equation}
\delta P \rho_s = \frac{\partial}{\partial y}\left(\frac{\rho_s(x_L)}{F'(x_L)}\delta F(x_L) - \frac{\rho_s(x_R)}{F'(x_R)} \delta 
F(x_R)\right)
 \label{e:cdf_pert3a}
\end{equation}

\noindent Combining equations (\ref{e:FP_bkw}) and (\ref{e:cdf_pert3a}), an expression for $\delta 
\bar{J}$ in terms of a mapping function perturbation $\delta F$ is obtained:

\begin{equation}
\delta \bar{J} = \int _0^1 \phi(y) \frac{\partial}{\partial y} 
\left(\frac{\rho_s(x_L)}{F'(x_L)}\delta F(x_L) - \frac{\rho_s(x_R)}{F'(x_R)} \delta 
F(x_R)\right)dy
\label{e:Jbar2}
\end{equation}

\noindent If $F$ and $\delta F$ are symmetric, $F'(x_L)$ will be positive and $F'(x_R)$ 
will be negative, therefore (\ref{e:density_mapping}) can be rewritten as:

\begin{equation}
\rho_s(F(x)) = \frac{1}{F'(x_L)}\rho_s(x_L) - \frac{1}{F'(x_R)}\rho_s(x_R)
\label{e:density_mapping2}
\end{equation}

\noindent Combining equations (\ref{e:Jbar2}) and (\ref{e:density_mapping2}):

\begin{equation}
\delta \bar{J} = \int _0^1 \phi(y) \frac{\partial}{\partial y}( \rho_s(y)\delta 
F(F^{-1}(y)))dy
\label{e:Jbar3}
\end{equation}

This is consistent with the equation for the density derivative in 
\cite{Rulle:1997:SRB}. To compute the gradient with respect to some parameter $\xi$, 
substitute $\frac{\partial F}{\partial \xi} \delta \xi$ for $\delta F$ in equation (\ref{e:Jbar3}) 
and divide through by $\delta \xi$:

\begin{equation}
\frac{\partial \bar{J}}{\partial \xi} = \lim_{\delta \xi \to 0} \frac{\delta \bar{J}}{\delta \xi} = \int 
_0^1 \phi(y) \frac{\partial}{\partial y} 
\left(\frac{\rho_s(x_L)}{F'(x_L)}\frac{\partial F}{\partial\xi}|_{x_L} - 
\frac{\rho_s(x_R)}{F'(x_R)} \frac{\partial F}{\partial\xi}|_{x_R}\right)dy
\label{e:JGrad}
\end{equation}

\noindent Or for the symmetric case:

\begin{equation}
\frac{\partial \bar{J}}{\partial \xi} = \int _0^1 \phi(y) 
\frac{\partial}{\partial y}\left(\rho_s(y)\frac{\partial 
F}{\partial\xi}|_{F^{-1}(y)}\right)dy
\label{e:JGrad_sym}
\end{equation}

Finally, care needs to be taken when discretizing the density adjoint equations.  The first eigenvalue of the discrete operator $P_n$ is not exactly
one and can change when the system is perturbed.  Because of this, an additional adjoint 
equation is required for $\lambda$ to compute the discrete density adjoint (see appendix \ref{ap:1Ddisc} for a derivation):

\begin{equation}
\left[ \begin{array}{cc} P_n^T - \lambda I & -\hv \\ -\hr_s^T & 0  \end{array} \right] \left[ 
\begin{array}{c} \hp \\ \eta  \end{array} \right] = \left[ \begin{array}{c} \underline{J} 
\\ 0 \end{array} \right]
\label{e:adj}
\end{equation}

\noindent Where $\eta$ is the adjoint of $\lambda$ and can be shown to
be equal to $\overline{J}$ in the continuous limit.  

\subsection{Algorithm Summary}

To compute some gradient $\frac{\partial \overline{J}}{\partial \xi}$, the 
following algorithm was used:

\begin{enumerate}
\item Compute the inverse of the mapping function $F(x)$ using Newton's Method.
\item Construct the matrix $P_n$ using the equations outlined in 
section \ref{ss:stat_dist}.  
\item Determine the stationary density $\hr_s$ using a power method.  Also 
determine the left eigenvector $v$ corresponding to the eigenvalue $\lambda$ of 
$\hr_s$.  
\item Compute the adjoint variable $\hp$ by solving (\ref{e:adj}).  To solve 
(\ref{e:adj}), be sure to take advantage of the sparseness of $P_n$.  
\item Compute the gradient using (\ref{e:JGrad}) or (\ref{e:JGrad_sym}).  
Approximate the $y$-derivative with a 2nd order center finite difference scheme.  
\end{enumerate}

\subsection{Density adjoint for the cusp map}
\label{ss:cusp_results}

The sensitivity of the mean $\bar{x}$ with respect 
to the parameter $\xi$ for the cusp map was computed.  For comparison, 
$\frac{\partial \bar{x}}{\partial \xi}$ was also computed using 1st order 
finite differences of equation (\ref{e:phase_space_avg}) adapted for the 1D case and 
discretized in n nodes:

\begin{equation}
\bar{x} = (1/n)x^T \hr_s
\end{equation}

$\overline{x}(\xi)$ was found to be sufficiently smooth to ensure accurate gradient computations using finite differences.  However, it is important to note this is not always the case \cite{Lea:2000:climate_sens}.

\begin{figure}[htb!] \centering
\includegraphics[width=3.2in]{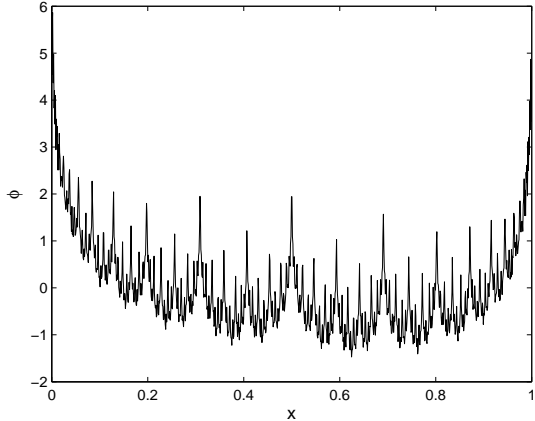}
\caption{Adjoint density $\phi$ for the cusp map with $\xi = 0.5$, generated using 1024 
nodes.  }
\label{f:cusp_adj}
\end{figure}

Figure \ref{f:cusp_adj} shows the density adjoint distribution for the cusp map with 
$\xi =0.5$.  The density adjoint is almost discontinuous, so small perturbations to stationary density $\rho_s$ can 
have large effects on the objective function.  Interestingly there is a fractal structure to the density adjoint.  This 
arises from the adjoint being computed backwards in time with the 
operator $P_n^T$.  $P_n$ folds and stretches density distributions, so $P_n^T$ duplicates 
and compresses features of adjoint density distributions.  This fractal structure arises
 because of the cusp map's ``peak'' at $x = 0.5$, which causes the folding and stretching.  



\begin{figure}[htbp!]
  \centering
  \begin{minipage}[b]{3.125in}
    \includegraphics[width=3.2in]{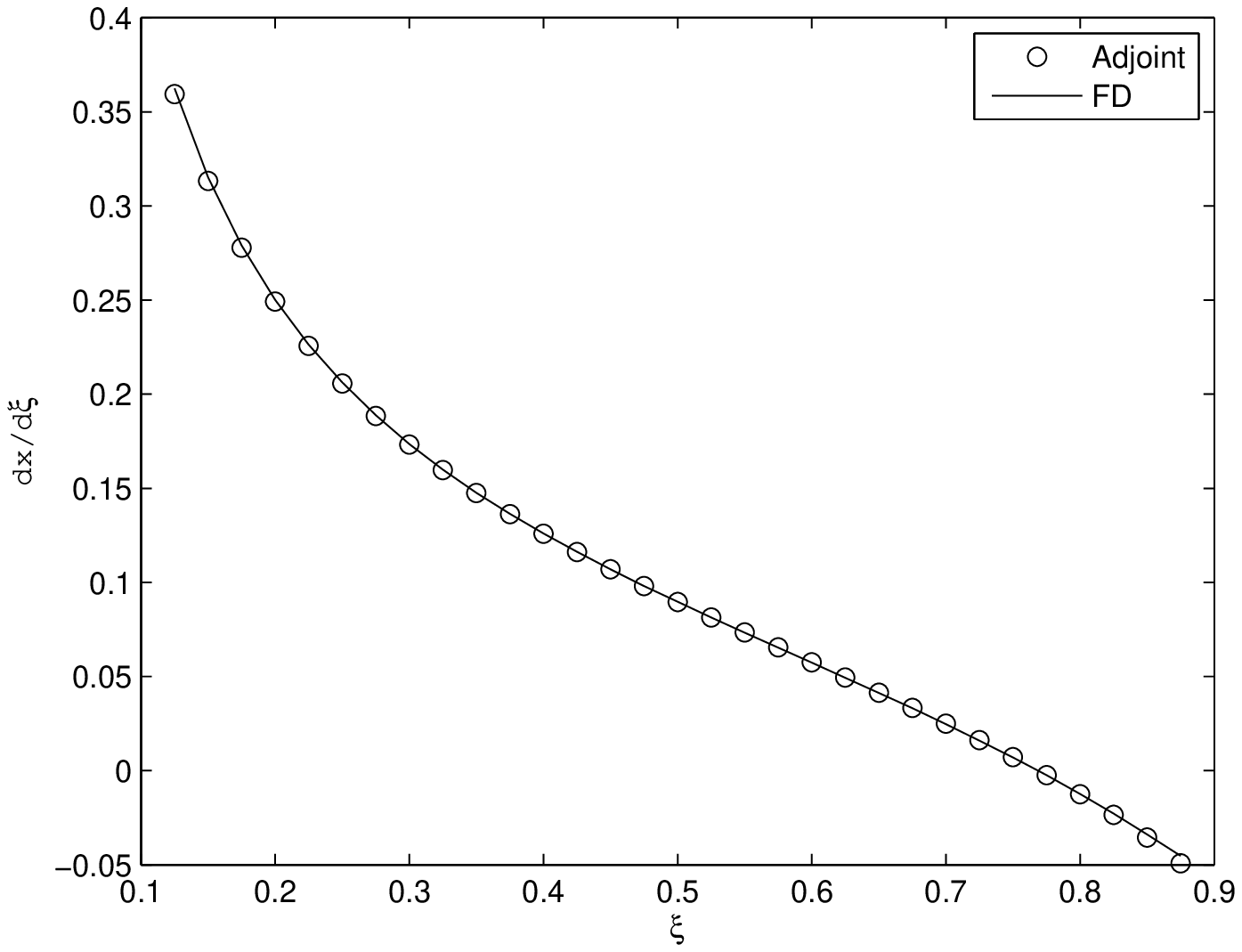}
    \caption{Comparison of gradients computed using the adjoint method and the 
finite difference method.  1D space between 0 and 1 was discretized using 256 
nodes.  }
    \label{f:adj_compare}
  \end{minipage}
\hspace{0.125 in}
  \begin{minipage}[b]{3.125in}
    \includegraphics[width = 3.2in]{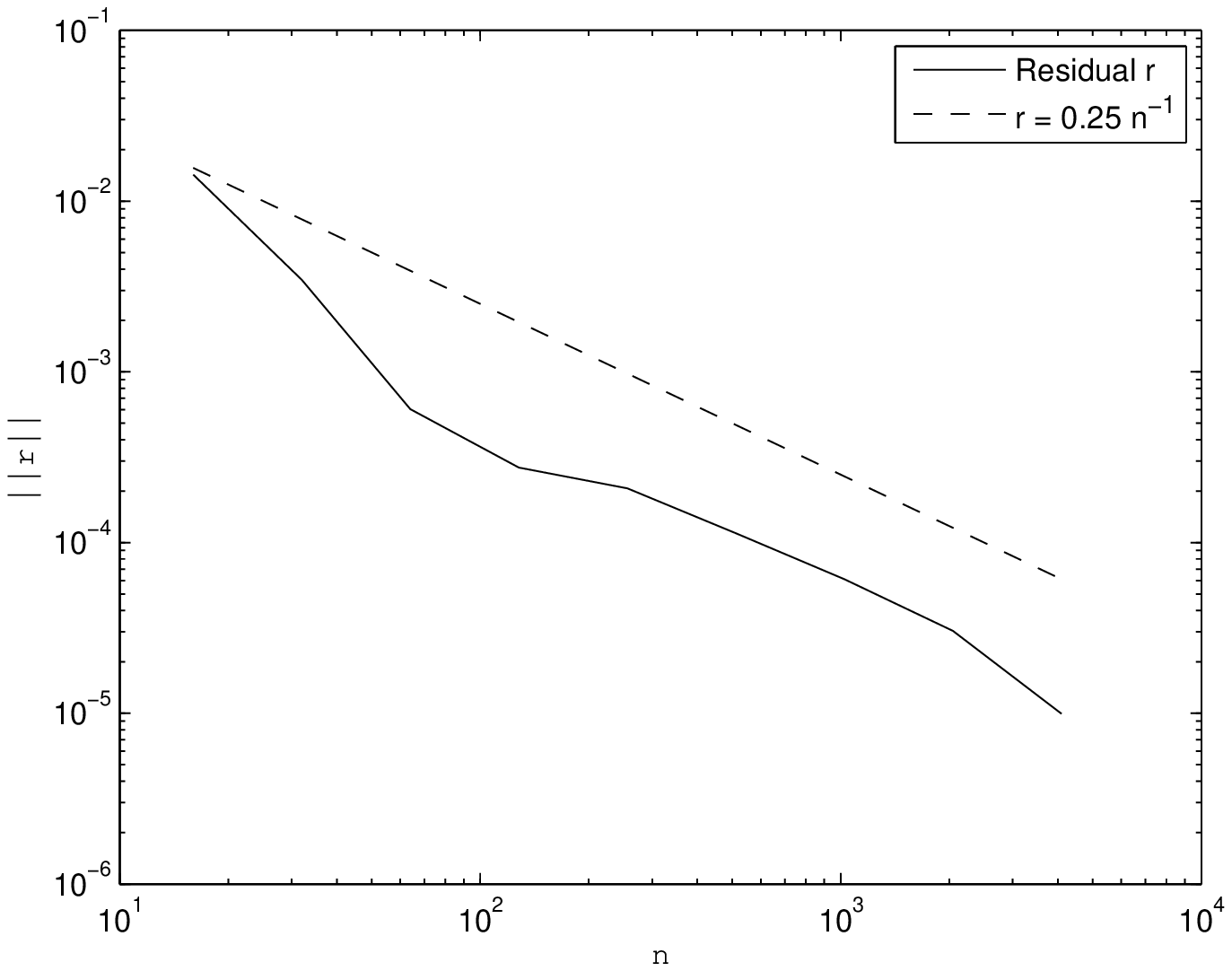}
\caption{Convergence of the residual $r$ of $\frac{\partial 
\overline{J}}{\partial \xi}$ with the number of nodes n for $\xi=0.5$.  The 
residual was calculated by taking the L2 norm of the difference between 
the gradient for $n$ nodes and 8096 nodes.  }
\label{f:1Dconvergence}
  \end{minipage}
\end{figure}


Despite the additional numerical dissipation, the 1D density adjoint computes
accurate gradient values.  Figure \ref{f:adj_compare} shows the adjoint and finite difference 
 computed gradients match up well visually.  It was found that the adjoint method predicts 
the gradient within 5\% of the finite difference calculation for most values of $\xi$. The order of convergence of the gradient varied 
slightly with $\xi$ and was typically around 1.15, as in figure \ref{f:1Dconvergence}.  




\section{Density adjoint for continuous chaos}
\label{s:densitylorenz}
The following section uses the Lorenz system as an example to illustrate the density adjoint method.  The method was used to compute the sensitivity 
of $\overline{z}$ to the parameters $s$, $r$, $b$ and $z_0$ in the 
Lorenz system:

\begin{align*}
\dot{x} &= s (y-x) \\
\dot{y} &= -x(z-z_0) + rx - y \\
\dot{z} &= xy - b(z-z_0)
\end{align*}

The parameters were set to their canonical values of $s = 10$, $r = 28$, 
$b=8/3$ and $z_0 = 0$. The Lorenz attractor with these parameters has a fractal 
dimension of roughly 2.05, so the attractor was approximated as a 2D surface in 3D 
phase space.  

\subsection{Computing Stationary Density}
\label{ss:stat_dist2}
In multiple dimensions, one could build a discrete Frobenius-Perron operator $P_n$ as in 
the 1D case.  $P_n$ would be an $\mathcal{M}$ by $\mathcal{M}$ matrix, where $\mathcal{M}$ is number of cells or nodes 
used to discretize the strange attractor.  To reduce the size of the matrix $P_n$, the matrix is built 
for a Poincar\'e section.  In this case, $P_n$ is $M$ by $M$, where $M$ is the number of nodes in the
Poincar\'e section, which is typically a small fraction of the total number of nodes $\mathcal{M}$.    
For the Lorenz attractor, a good choice for the 
Poincar\'e section is a constant $z$ plane including the two non-zero unstable fixed 
points at $(\pm\sqrt{b(r-1)},\pm\sqrt{b(r-1)},r+z_0-1)$.  

\begin{figure}[htb!]
  \centering
  \begin{minipage}[b]{3.125in}
    \includegraphics[width=3.2in]{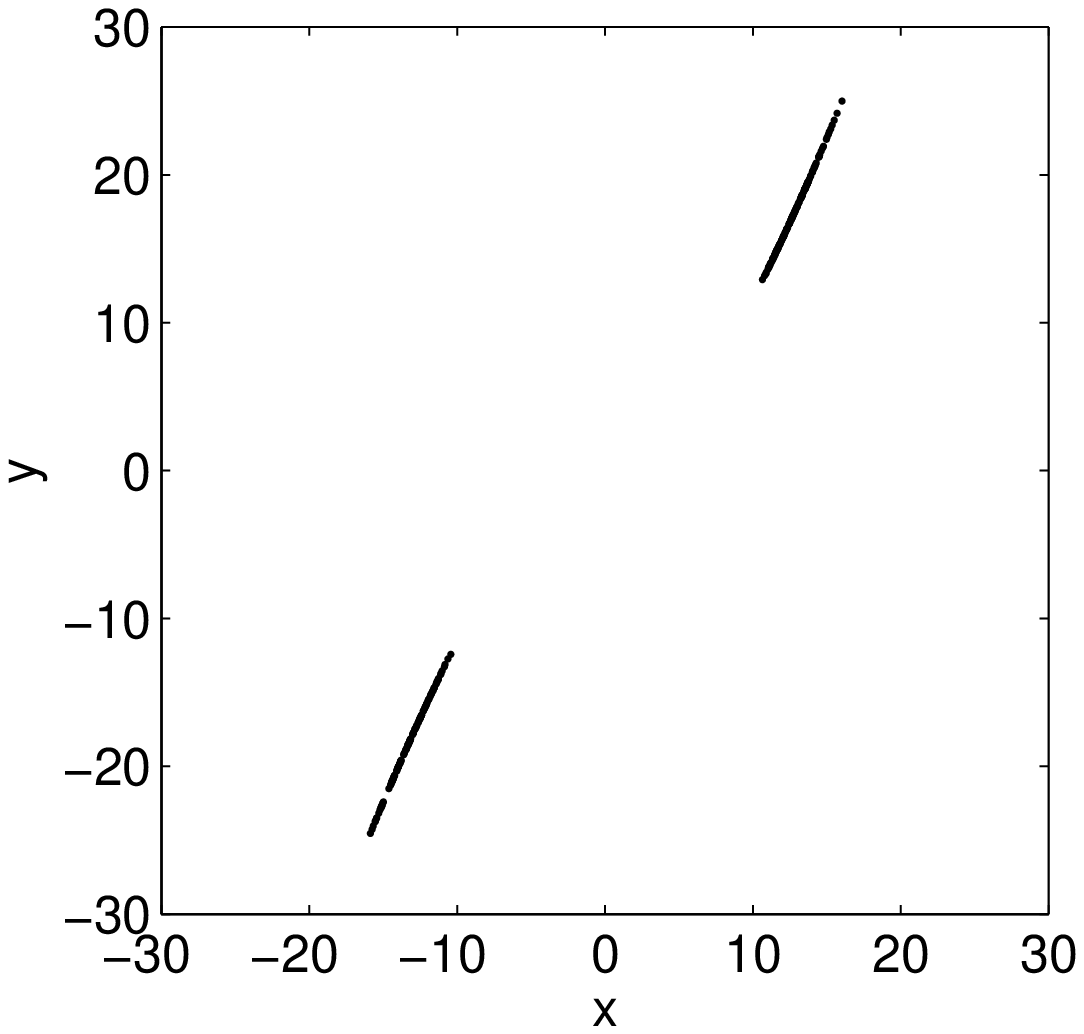}
    \caption{Poincar\'e Section at $z = 27$ for Lorenz attractor trajectories 
with $\frac{\partial z}{\partial t} > 0$.  }
    \label{f:lorenzPoincare2D}
  \end{minipage}
\hspace{0.125 in}
  \begin{minipage}[b]{3.125in}
    \includegraphics[width=3.2in]{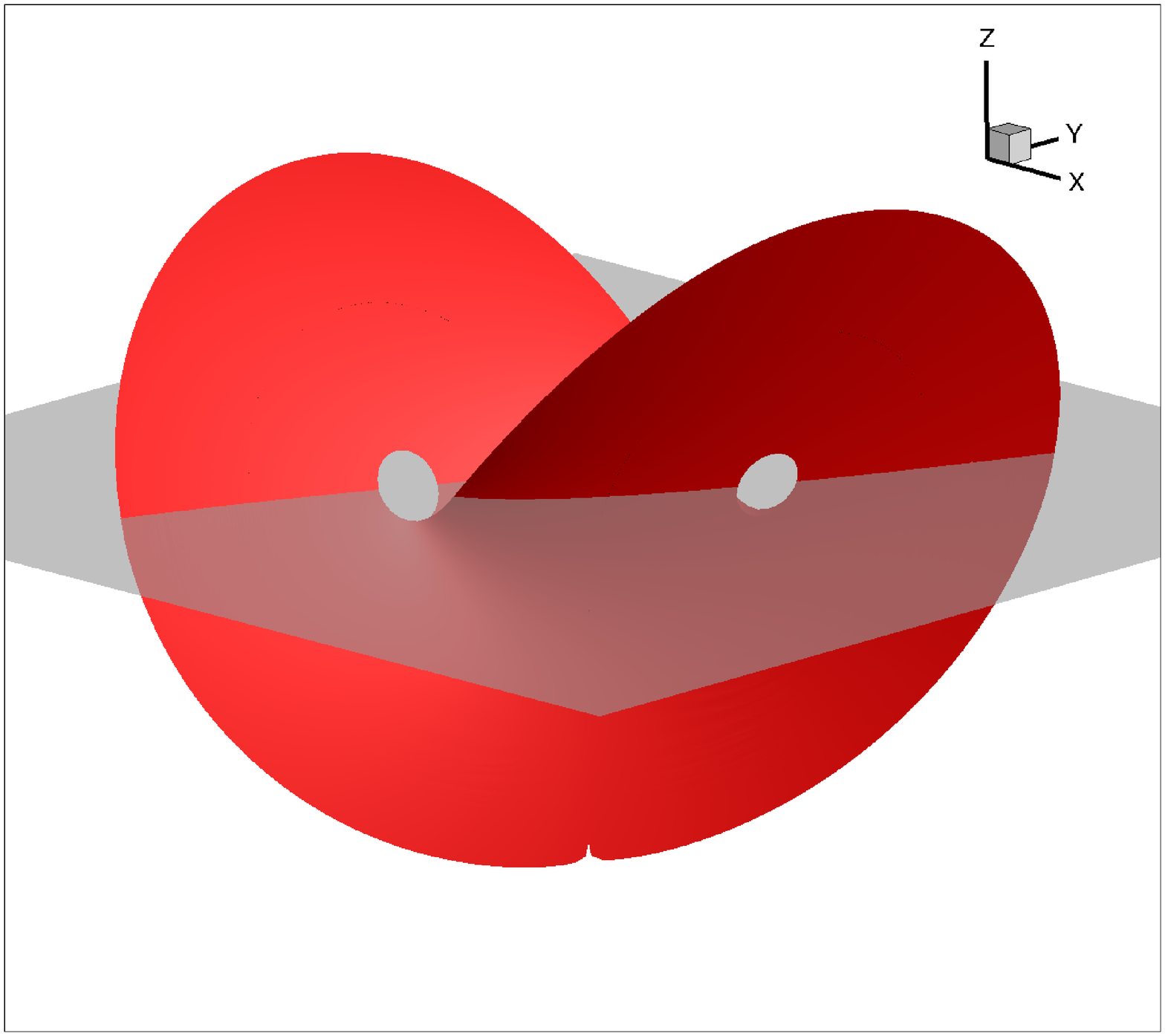}
    \caption{Three dimensional view of the 2D surface approximating the Lorenz Attractor and the Poincar\'e 
section at $z = 27$. }
    \label{f:lorenzPoincare3D}
  \end{minipage}
\end{figure}

This Poincar\'e section has an attractor cross-section that can be well approximated as
a 1D function of either $x$ or $y$.  Therefore solving for the density distribution at 
the Poincar\'e section is a 1D map problem and the stationary density 
distribution in the Poincar\'e section $\rho_0$ can be computed from the operator $P_n$ as shown in 
section \ref{s:density1d}.  For the Lorenz attractor, the starting positions of 
the streamlines, defined as the trajectories $\vec{x}(t)$ used to discretize the attractor, were determined using a 7th order polynomial curve fit through 
a Poincar\'e section taken from a trajectory with length $T = 10000$ time units.  

As the Poincar\'e section can be modeled with a polynomial curve fit, the
 attractor itself can be approximately modeled as a 2D surface, as shown in 
figure \ref{f:lorenzPoincare3D}.  

As the Lorenz attractor lies in a three dimensional phase space, vector notation is used in this section.  
A lower case symbol is a scalar (i.e. stationary density $\rho_s$), a symbol with 
an arrow overhead is a column vector (i.e. phase space position $\vec{x}$) and
a matrix/tensor is indicated by bold script (i.e. a Jacobian $\textbf{J}$).  

Unlike the 1D case, an explicit form of the mapping function is not available. 
Instead, a probability mass conservation equation is derived from the normalization
 and additivity axioms of probability.  The probability mass conservation equation 
can be used to compute the ratio between densities for a given "mapping", which can 
be used instead of the mapping function slope in the density mapping equation.  

A very helpful physical analogy to the conservation of probability on the 
attractor surface is the conservation of mass in a fluid flow.  Like mass, 
probability cannot be created or destroyed according to the nonnegativity 
and normalization axioms.  Therefore, the following equation holds:

\begin{equation}
\vec{\nabla}_s \cdot (\rho_s(\vec{x})\vec{f}(\vec{x})) = 0
\label{e:dens_cons}
\end{equation}

\noindent Where the gradient operator $\vec{\nabla}_s$ is an operator on the attractor manifold, 
$\vec{x}$ is a point in phase space $R^n$ and:

\[
\frac{d\vec{x}}{dt} = \vec{f}(\vec{x})
\]

\noindent is the system of equations governing the dynamical system of interest (i.e. the Lorenz system).  The 
physical analog of $\vec{f}(\vec{x})$ is a velocity field in a fluid flow, hence the name ``streamline'' for a phase space trajectory on the attractor.  $\vec{x}(t)$ gives the path that density $\rho$ ``flows'' on the attractor manifold.  

From (\ref{e:dens_cons}), a partial differential equation (PDE) governing the 
density distribution on an attractor can be derived.  From the chain rule\footnote{$\rho(\vec{x})$ can be any probability distribution on the attractor, including the stationary distribution $\rho_s(\vec{x}) = \lim_{t \to \infty} \rho(\vec{x})$}:
\[
\vec{f}(\vec{x}) \cdot \vec{\nabla}_s \rho(\vec{x},t) + \rho(\vec{x},t) \vec{\nabla}_s \cdot 
\vec{f}(\vec{x})= 0 
\]

The Lorenz attractor is approximated as a 2D surface, so two natural 
coordinates are used; $l$, which is in the direction of the ``velocity field'' 
defined by $\vec{f}(\vec{x})$ and $s$, which is orthogonal to $l$ but tangent to the 
attractor surface.  $l$ and $s$ will be referred to as the streamwise and spanwise 
directions respectively. $\hat{l}$ and $\hat{s}$ are unit vectors in the streamwise and spanwise
directions. Using thes definitions, the density PDE can be simplified:
\begin{equation}
|\vec{f}(\vec{x})| \frac{\partial \rho}{\partial l}= - \rho(\vec{x},t) \vec{\nabla}_s \cdot 
\vec{f}(\vec{x}) \quad \Rightarrow \quad 
\frac{\partial \rho}{\partial l}= - \frac{\rho(\vec{x},t)}{|\vec{f}(\vec{x})|} \vec{\nabla}_s \cdot \vec{f}(\vec{x})
\label{e:density_pde}
\end{equation}

\noindent Additionally, it can be shown that (see \ref{ap:div}):

\begin{equation}
 \vec{\nabla}_s  \cdot \vec{f}(\vec{x}) = \hat{l}^T \textbf{J}\hat{l} + \hat{s}^T \textbf{J} \hat{s}
 \label{e:attractor_div}
\end{equation}

\noindent Where $\textbf{J}$ is the Jacobian of $\vec{f}(\vec{x})$. Equation \eqref{e:attractor_div} can be substituted into equation \eqref{e:density_pde} to obtain:

\begin{equation}
\frac{\partial \rho}{\partial l}= - \frac{\rho(\vec{x},t)}{|\vec{f}(\vec{x})|} (\hat{l}^T \textbf{J}\hat{l} + \hat{s}^T \textbf{J} \hat{s})
\label{e:dens_PDEl}
\end{equation}

As $\hat{l}$ and $\vec{f}(\vec{x})$ are the same direction, and $d\vec{x}/dt=\vec{f}(\vec{x})$, the following relation holds for a  streamline $\vec{x}(t)$ on the attractor:

\[
\frac{dl}{dt} = |\vec{f}(\vec{x})|
\] 

\noindent Therefore,

\begin{equation}
\frac{\partial \rho}{\partial t}= - \rho(\vec{x},t) (\hat{l}^T \textbf{J}\hat{l} + \hat{s}^T \textbf{J} \hat{s})
\label{e:dens_PDE}
\end{equation}

As $\rho$ is invariant when multiplied by a constant, this equation can be 
time integrated along some streamline $\vec{x}(t)$ to find the ratio between density at different points on a given 
Poincar\'e section: 

\begin{equation}
 \log \frac{\rho(T)}{\rho_0} = -\int_0^T (\hat{l}^T \textbf{J}\hat{l} + \hat{s}^T \textbf{J} \hat{s})dt
\label{e:dens_ratio}
\end{equation}

\noindent Where $\rho_0$ is defined as the density at the beginning of the streamline starting at $\vec{x}_0 = \vec{x}(0)$ and  $\rho(T)$ is the density at $\vec{x}(T)$, where the streamline $\vec{x}(t)$ returns to the Poincar\'e section.  
Equation (\ref{e:dens_ratio}) is numerically integrated to find the 
ratio between the density at the beginning and end of $M$ streamlines.  
These ratios, along with the start and end positions in phase space of each streamline can be 
used to form a Frobenius-Perron operator $P_n$ with a first eigenvector corresponding to 
the stationary density distribution at the Poincar\'e section.  As the starting 
and ending positions of the streamlines will rarely match (i.e. $\vec{x}(T)_i \neq \vec{x}(0)_j$), linear interpolation is used 
as in the 1D case to compute the density ``flow'' between the starting and ending positions. 

\begin{figure}[htb!]
  \centering 
    \includegraphics[width=3.2in]{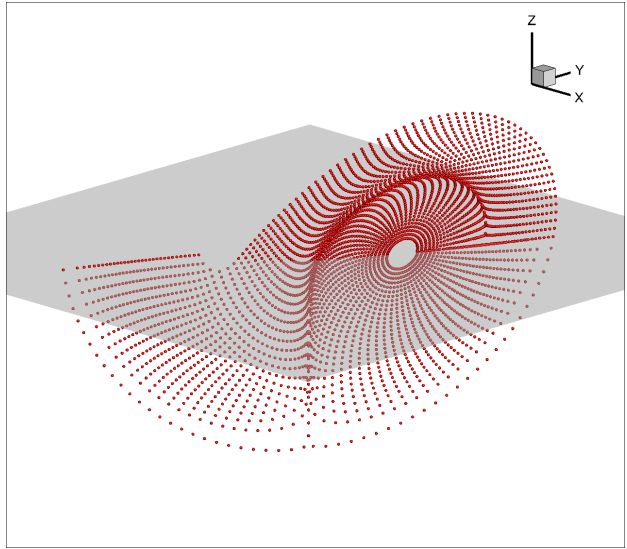}
    \caption{Node distribution corresponding to a 64 streamline by 64 
streamwise mesh for the Lorenz attractor. It was found that distributing 
the streamline starting positions so that there were more streamlines near the 
bifurcation increased the rate of convergence to the true density distribution. 
 }
    \label{f:lorenz64mesh}
\end{figure}
\begin{figure}[htb!]
  \centering
  \begin{minipage}[b]{3.125in}
    \includegraphics[width=3.2in]{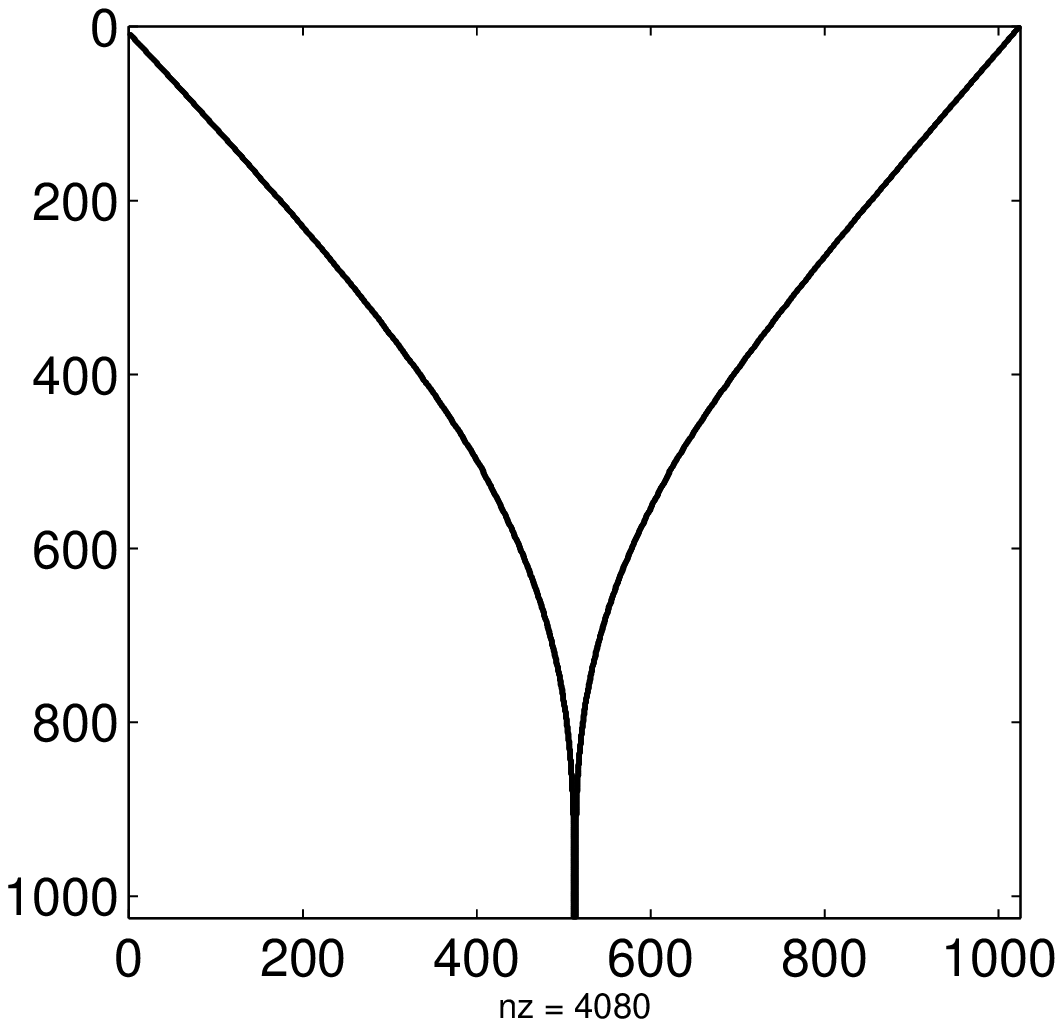}
    \caption{Transition Matrix $P_n$ structure for a roughly uniform streamline 
distribution.  Note the similarity of this matrix to that for the Cusp map}
    \label{f:lorenz_matrix}
  \end{minipage}
\hspace{0.125 in}
  \begin{minipage}[b]{3.125in}
  	\includegraphics[width=3.2in]{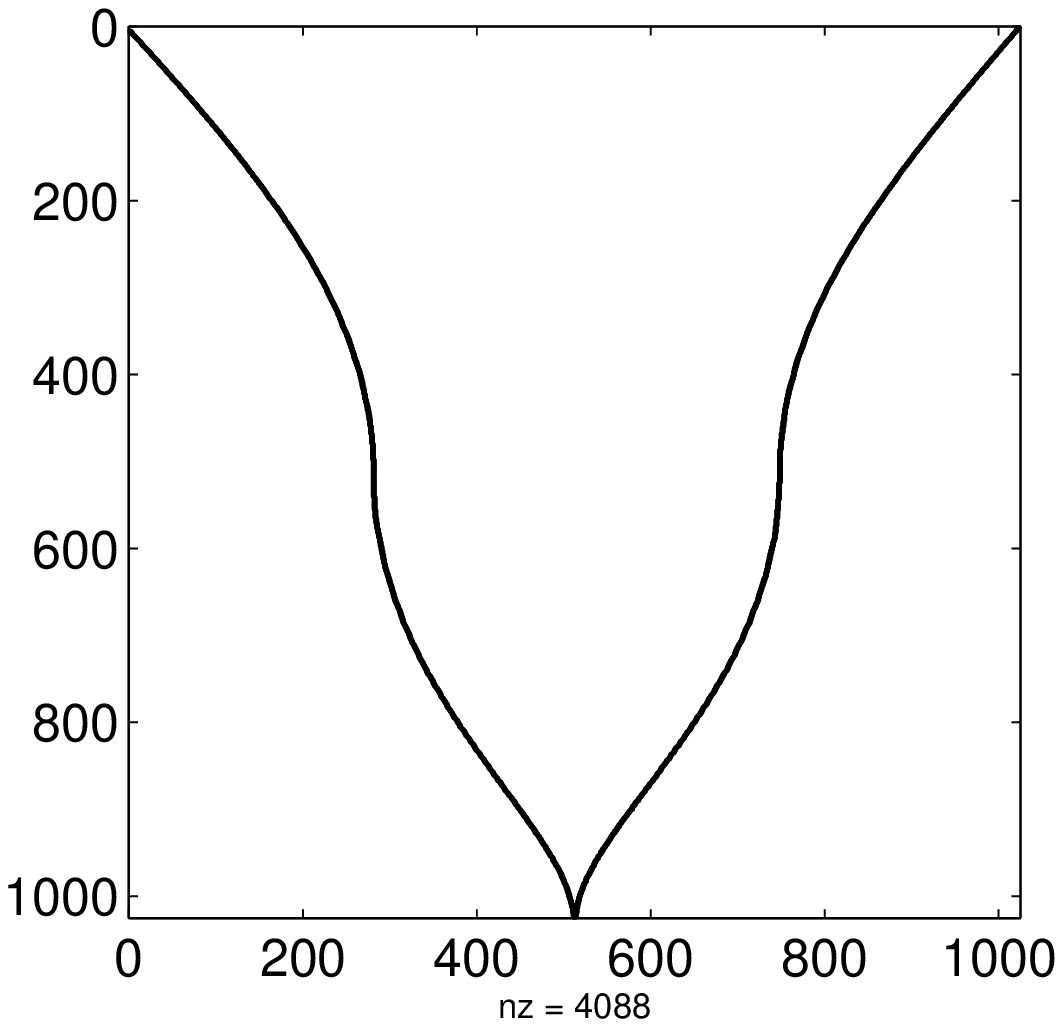}
    \caption{Transition Matrix $P_n$ structure for a non-uniform streamline 
distribution with more streamlines starting near $x=13,y=18,z=27$.  }
    \label{f:lorenz_matrix_NU}
  \end{minipage}

    \begin{minipage}[b]{3.125in}
	\includegraphics[width=3.2in]{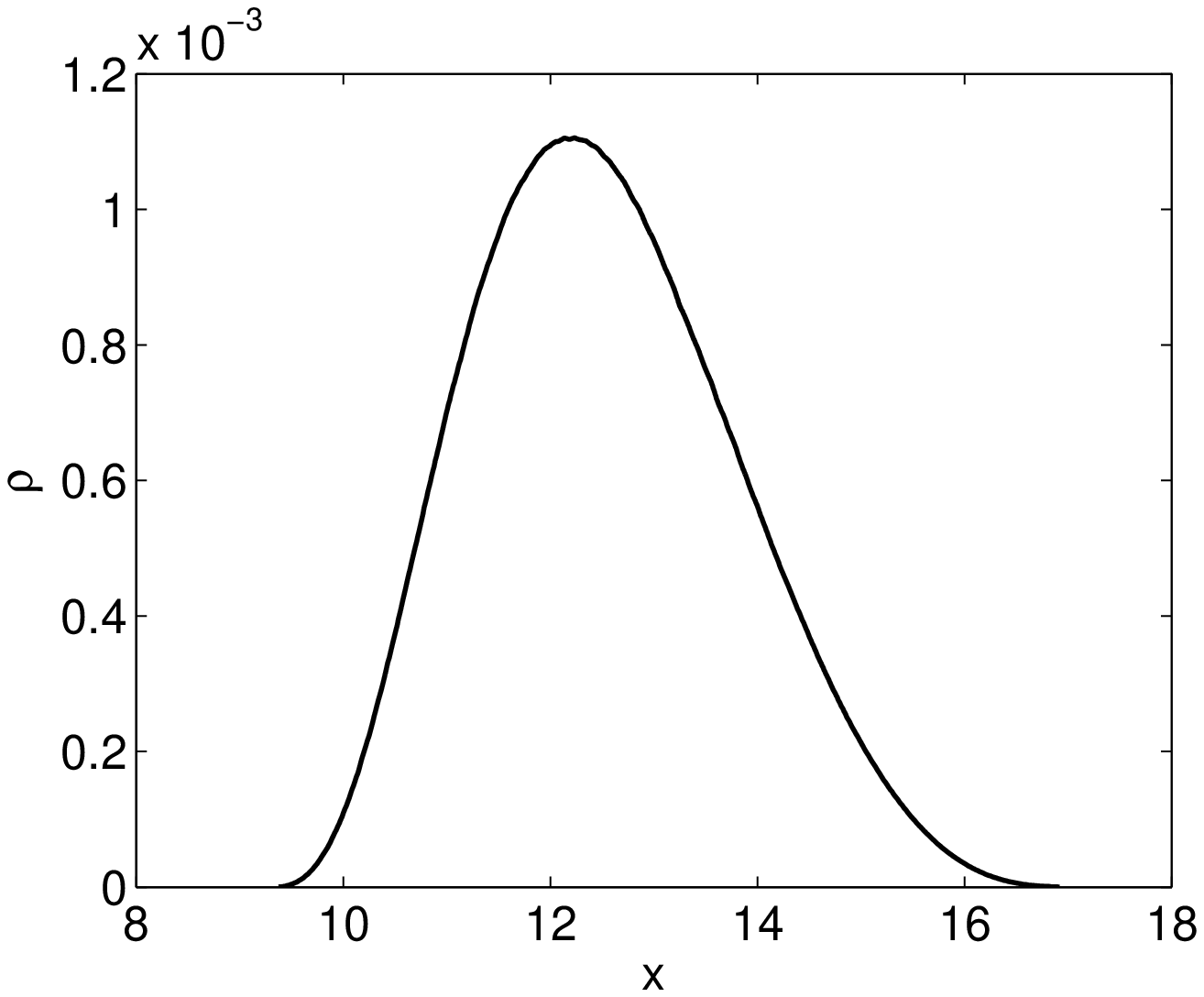}
    \caption{Density $\rho_s$ versus $y$ on the Poincare Section at $z = 27$. 512 streamlines were used to form $P_n$.  }
    \label{f:poincare_rho}
	\end{minipage} 
  \hspace{0.125 in} 
    \begin{minipage}[b]{3.125in}
    \includegraphics[width=3.2in]{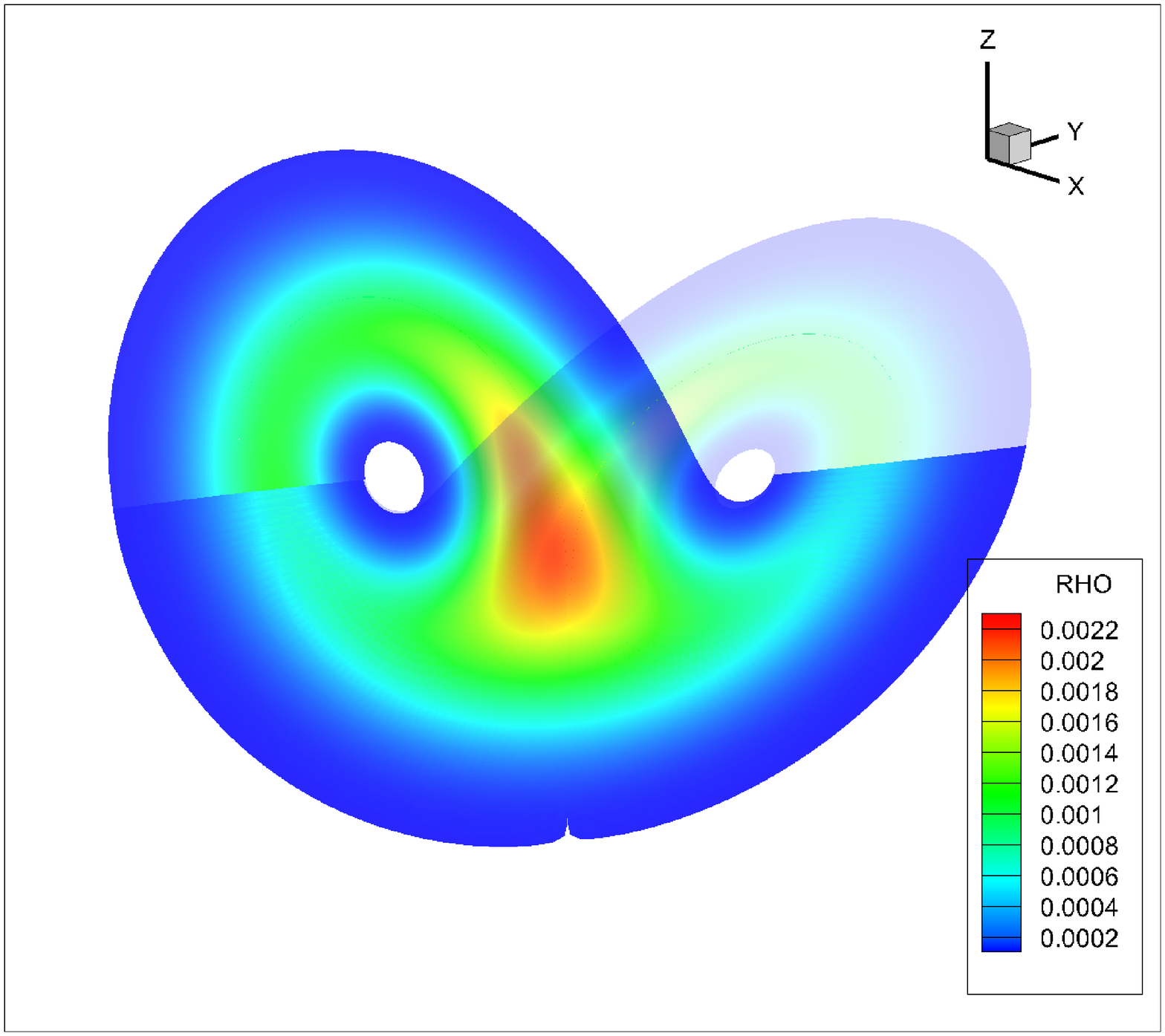}
    \caption{Density distribution on the surface of the Lorenz attractor for a 
512 by 128 mesh.  }
    \label{f:lorenz_rho}
    \end{minipage} 
    
\end{figure}

By the symmetry of the Lorenz system, the Poincar\'e plane intersections for 
the Lorenz attractor are 180 degree rotational translations of one another as 
is evident in figure \ref{f:lorenzPoincare2D}, where it can be seen that $x = -x$ and 
$y = -y$.  
This symmetry of the attractor can be exploited for 
lower computational costs.  If the attractor is discretized with streamlines 
starting along the Poincar\'e section in 
the first quadrant ($x>0, \ y>0$), a portion of the streamlines return to the first quadrant 
and a portion go to the third quadrant ($x<0, \ y<0$), as seen in figure \ref{f:lorenz64mesh}. 
 By symmetry, the streamlines running from the first to the third quadrant are 
the same as those from the third to the first rotated 180 degrees about the z-axis.
  This means that the density flux from the third quadrant is the same as the 
density flux to the third quadrant. The density flow from returning and incoming 
streamlines make up two sides of the transition matrix $P_n$, as shown in figures 
\ref{f:lorenz_matrix} and \ref{f:lorenz_matrix_NU}.  


Because the linear interpolation scheme does not ensure conservation of 
probability mass, the Poincar\'e stationary distribution $\rho_0(\vec{x}_0)$ as computed 
using a power method is not properly normalized. This is because the first 
eigenvalue is not equal to one as it would be if probability mass was conserved.  To 
normalize $\rho_0(\vec{x}_0)$ begin with the density over the entire attractor, $\rho(\vec{x})$ and use the fact that
 $dl = |\vec{f}(\vec{x})|dt$:

\begin{equation}
\iint \rho \ dl ds  = \iint \rho |\vec{f}(\vec{x})| dt ds = 1
\label{e:dens_int}
\end{equation}

\noindent Conservation of probability mass along a streamline can be written as:

\[
\rho |\vec{f}(\vec{x})| ds = \rho_0 |\vec{f}(\vec{x}_0) \times \hat{s}_0| ds_0
\]

Where $ds$ is the width of the streamline at a given $\vec{x}$, $ds_0$ is its 
initial width (width at $\vec{x}_0$), $\vec{f}(\vec{x}_0)$ is the "initial velocity" and $\hat{s}_0$ is the 
initial spanwise direction.  Substituting into equation (\ref{e:dens_int}):

\begin{align*}
\iint_0^T \rho_0 |\vec{f}(\vec{x}_0) \times \hat{s}_0| ds_0 dt  &= 1 \\
\int  \rho_0 \left(\int_0^T dt \right) |\vec{f}(\vec{x}_0) \times \hat{s}_0| ds_0  &= 1 \\
\int \rho_0 T |\vec{f}(\vec{x}_0) \times \hat{s}_0| ds_0  &= 1
\end{align*}

\noindent The discretized form of this equation, which can be used to normalize $\rho_0$, is:

\[
\hr_0^T \hv = 1
\]

\noindent Where:

\begin{equation}
\hv_i = T |f(x_0) \times \hat{s}_0| ds_0
\label{e:lead_eig_v}
\end{equation}

\noindent and $T$ is the total time a particle spends along streamline $i$. It can be shown 
that $v$ is the leading eigenvector of $P_n$ corresponding to $\lambda \approx 1$. 

Figure \ref{f:poincare_rho} shows the Poincar\'e section stationary distribution 
$\rho_0(\vec{x}_0)$ for the Lorenz attractor.  Like the stationary distribution for the cusp map it is 
smooth and continuous. 

Once the Poincar\'e stationary distribution $\rho_0(\vec{x}_0)$ is computed and 
normalized, the stationary density distribution over the entire attractor is 
computed by integrating equation (\ref{e:dens_PDE}) for each streamline with 
$\rho_0(\vec{x}_0)$ as the initial value of the stationary distribution $\rho_s$ along a streamline starting at 
$\vec{x}_0$.  

The density distribution computed for the Lorenz attractor is shown in figure 
\ref{f:lorenz_rho}.  The apparent discontinuity results from the intersection 
of the two branches of the attractor.  The sum of the density distribution on 
the intersection of these two branches is equal to the distribution on the 
Poincar\'e section.  

\begin{figure}[htbp]
  \centering
    \includegraphics[width=3.2in]{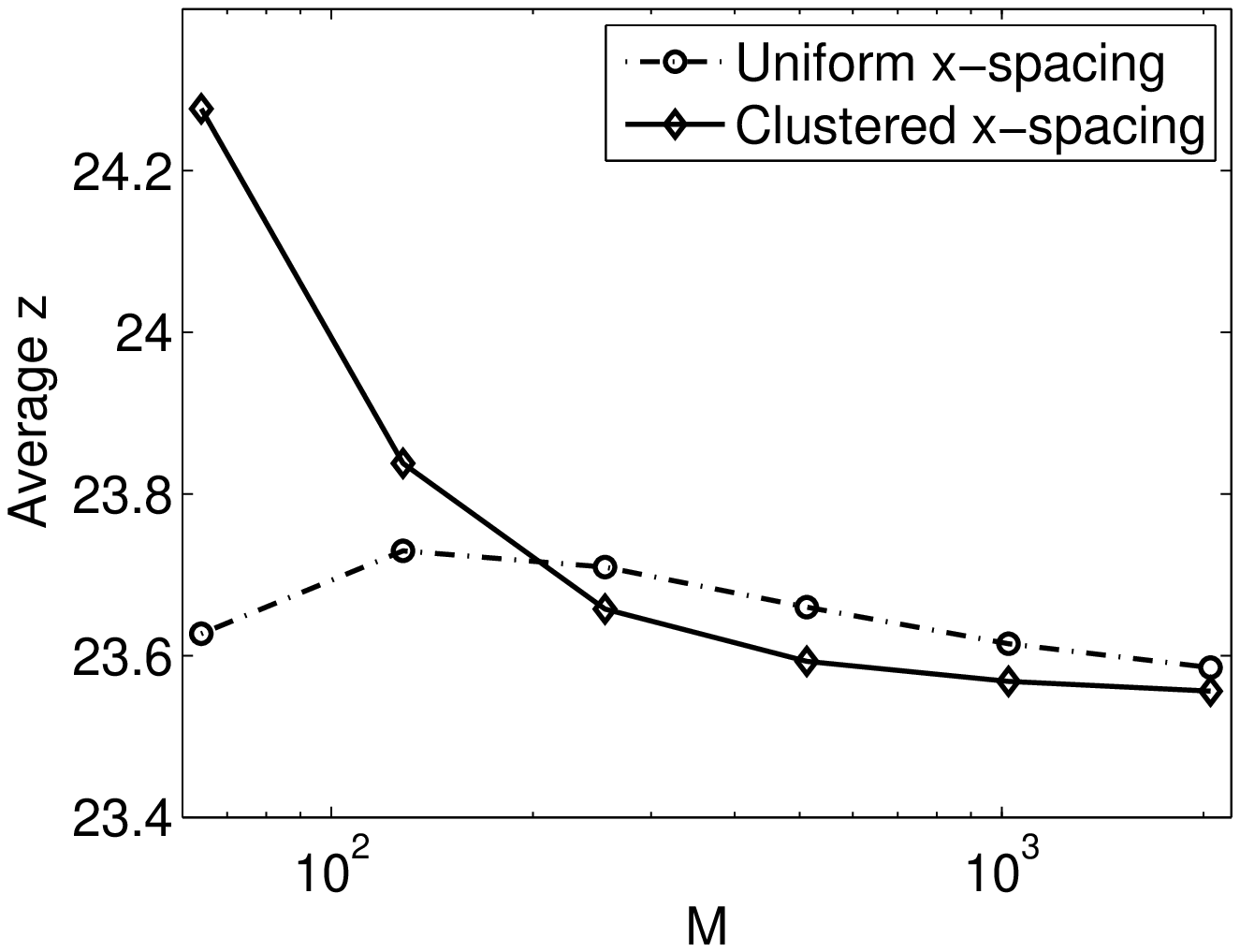}
    \caption{Convergence of $\overline{z}$ for two different streamline start 
position distributions, where $M$ is the number of streamlines.  The clustered 
distribution has streamlines clustered near the bifurcation of the attractor.  }
    \label{f:zbar_conv}

\end{figure}

Figure \ref{f:zbar_conv} shows that $\overline{z}=23.6$ 
from the density distribution, which is consistent with the value 
$\overline{z}=23.550$ found using ensemble averages of long phase space trajectories.

\subsection{Computing the Density Adjoint}

%

As the negative Lyapunov exponent for the Lorenz attractor has a large 
magnitude relative to the positive Lyapunov exponent, it can be assumed that perturbations 
to the long-time averaged quantity $\overline{J}$ arise mainly from perturbations 
to the stationary density $\delta \rho$ as opposed to perturbations to the 
attractor manifold:

%

\begin{equation}
 \delta\overline{J} = \iint J(\vec{x})\delta \rho \ dl ds
\label{e:delta_J}
\end{equation}

The adjoint density equation can be found using equation \eqref{e:delta_J} and 
the linearization of equation \eqref{e:dens_cons} (see \ref{ap:NDcont} for 
detailed derivation):

\begin{equation}
\frac{\partial \phi}{\partial t} = J(\vec{x}) -\overline{J} 
\label{e:lorenz_adj}
\end{equation}

\noindent Perturbations to $\overline{J}$ can then be computed using:

\begin{equation}
 \delta \overline{J} = \iint \phi \, \vec{\nabla}_s \cdot (\rho_s \delta \vec{f}) \ dl ds
 \label{e:adj_grad}
\end{equation}

\noindent Therefore gradients with respect to some parameter $\xi$ are:

\[
\frac{\partial \overline{J}}{\partial \xi} = \iint \phi \, \vec{\nabla}_s \cdot 
\left(\rho_s \frac{\partial \vec{f}}{\partial \xi}\right) dl ds
\]

\noindent To derive the adjoint equations for a numerical scheme,  first consider equation
(\ref{e:phase_space_avg}):

\[
\overline{J} = \iint J(\vec{x}) \rho(\vec{x}) dl ds
\]

\noindent This can be rewritten as:

\[
\overline{J} = \iint_0^T J(t) dt \rho_0 |\vec{f}(\vec{x}_0) \times \hat{s}_0| ds_0 
\]

\noindent Defining $\hJ_i = \int_0^T J(t) dt$ for streamline $i$ ``flowing'' from the Poincar\'e section, the above 
equation has the discretized form:

\begin{equation}
\overline{J} = \hJ^T D \hr_0
\label{e:avg_frm_rho0}
\end{equation}

Where $D$ is a diagonal matrix with $|\vec{f}(\vec{x}_0) \times \hat{s}_0| ds_0$ for the 
ith streamline along the main diagonal. Using $D$ to rescale $P_n$, the adjoint equation for 
$\hr_0$ can be derived for the 1D Poincar\'e map as in section \ref{s:density1d} 
(see \ref{ap:NDdisc} for a detailed derivation):

\begin{equation}
\left[ \begin{array}{cc} (D^{-1}P^T D-\lambda I) & -D^{-1} \hv \\ \hr_s^T D & 
0  \end{array} \right] \left[ \begin{array}{c} \hp_0 \\ \overline{J}  
\end{array} \right] = \left[ \begin{array}{c} \hJ \\ 0 \end{array} 
\right]
\label{e:disc_adj1}
\end{equation}

$\lambda$ is included in equation \eqref{e:disc_adj1} because it is not exactly one in 
practice.  The adjoint density along the Poincar\'e section $\phi_0(\vec{x}_0)$ is computed using 
\eqref{e:disc_adj1}.  Then equation \eqref{e:lorenz_adj} is integrated to compute the 
adjoint along each streamline, using $\phi_0(\vec{x}_0)$ as the initial value.

\noindent Gradients can be computed by discretizing equation (\ref{e:adj_grad}):

\begin{equation}
 \frac{\partial \overline{J}}{\partial \xi} \approx \sum_{k=0}^{N} \phi_k \, 
[\vec{\nabla}_s \cdot \left(\rho_s \frac{\partial \vec{f}}{\partial \xi}\right)]_k dA_k
\label{e:lorenz_grad}
\end{equation}

Where $dA_k$ is the attractor manifold ``area'' corresponding to the kth node.  This can be computed 
by integrating a differential equation formed using conservation of probability mass
(see \ref{ap:attractor_area}).  The quantity $[\vec{\nabla}_s \cdot \left(\rho_s \frac{\partial f}{\partial \xi}\right)]_k$
can be computed by finite differences as $\rho_s$ is known for each node and 
$\frac{\partial \vec{f}}{\partial \xi}$ can be found analytically for each node 
(see \ref{ap:NDgrad} for a detailed derivation). 

\subsection{Algorithm Summary}
\label{ss:alg_ND}
To compute some gradient $\frac{\partial \overline{J}}{\partial \xi}$, the 
following algorithm was used:

\begin{enumerate}
\item Find a Poincar\'e Section for the attractor such that the intersections 
trace an approximately one to one function, as seen in figure 
\ref{f:lorenzPoincare2D}.  Find a curve fit for these intersections.  
\item Construct the matrix $P_n$ with a loop, by integrating (\ref{e:dens_PDE}) 
along a set of streamlines originating and terminating at the Poincar\'e Plane 
from step 1.  
\item Determine the stationary density $\rho_0$ on the Poincar\'e plane using a 
power method. Smooth this distribution using a low-pass filter if necessary.  
\item Compute $\overline{J}$ using the following equation:
\[
\overline{J} = \rho_0^T D \mathcal{J}_s
\]   
\item Determine the left eigenvector $v$ corresponding to the eigenvalue 
$\lambda$ of $\rho_0$ using (\ref{e:lead_eig_v}).  
\item Compute the Poincar\'e Plane adjoint density $\phi_0$ by solving (\ref{e:adj}).  To solve (\ref{e:disc_adj1}), be sure to take advantage of the 
sparseness of $P_n$. 
\item Using $\rho_0$ and $\phi_0$ as initial values, integrate 
(\ref{e:dens_PDE}) and (\ref{e:lorenz_adj}) along each streamline to find 
$\rho_s$ and $\phi$ for the entire attractor.  
\item Find $\frac{\partial \vec{f}}{\partial \xi}$ analytically and calculate its 
value at all nodes.  
\item Compute the gradient using (\ref{e:lorenz_grad}).   
\end{enumerate}

\subsection{Density adjoint for the Lorenz system}

\begin{figure}[htbp]
  \centering
  \begin{minipage}[b]{3.15in}
    \includegraphics[width=3.2in]{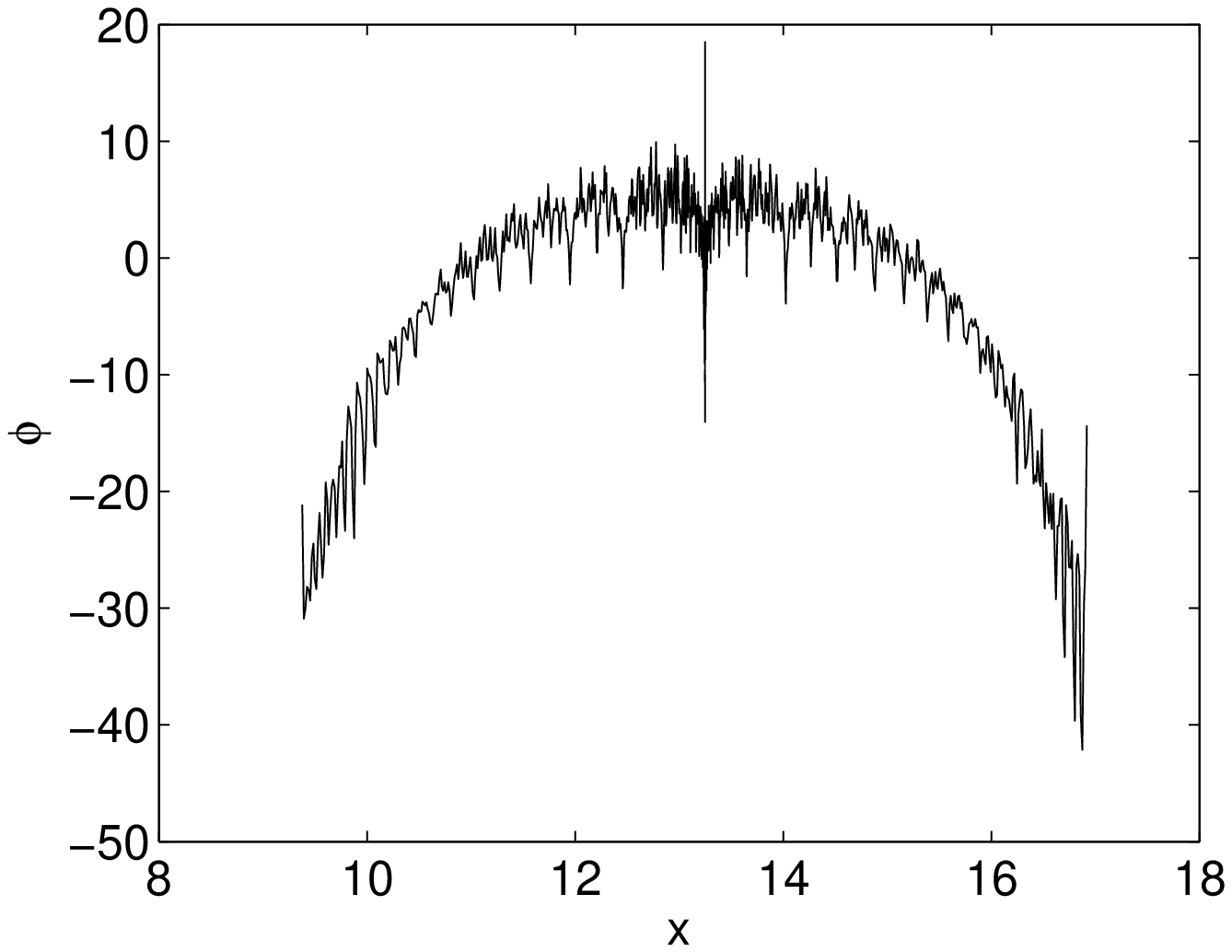}
    \caption{Adjoint $\phi$ versus $x$ on the Poincar\'e Section at $z = 27$. 1024 
streamlines were used to form $P_n$.  }
    \label{f:poincare_phi}
  \end{minipage}
\hspace{0.125 in}
  \begin{minipage}[b]{3.1in}
  
    \includegraphics[width=3.2in]{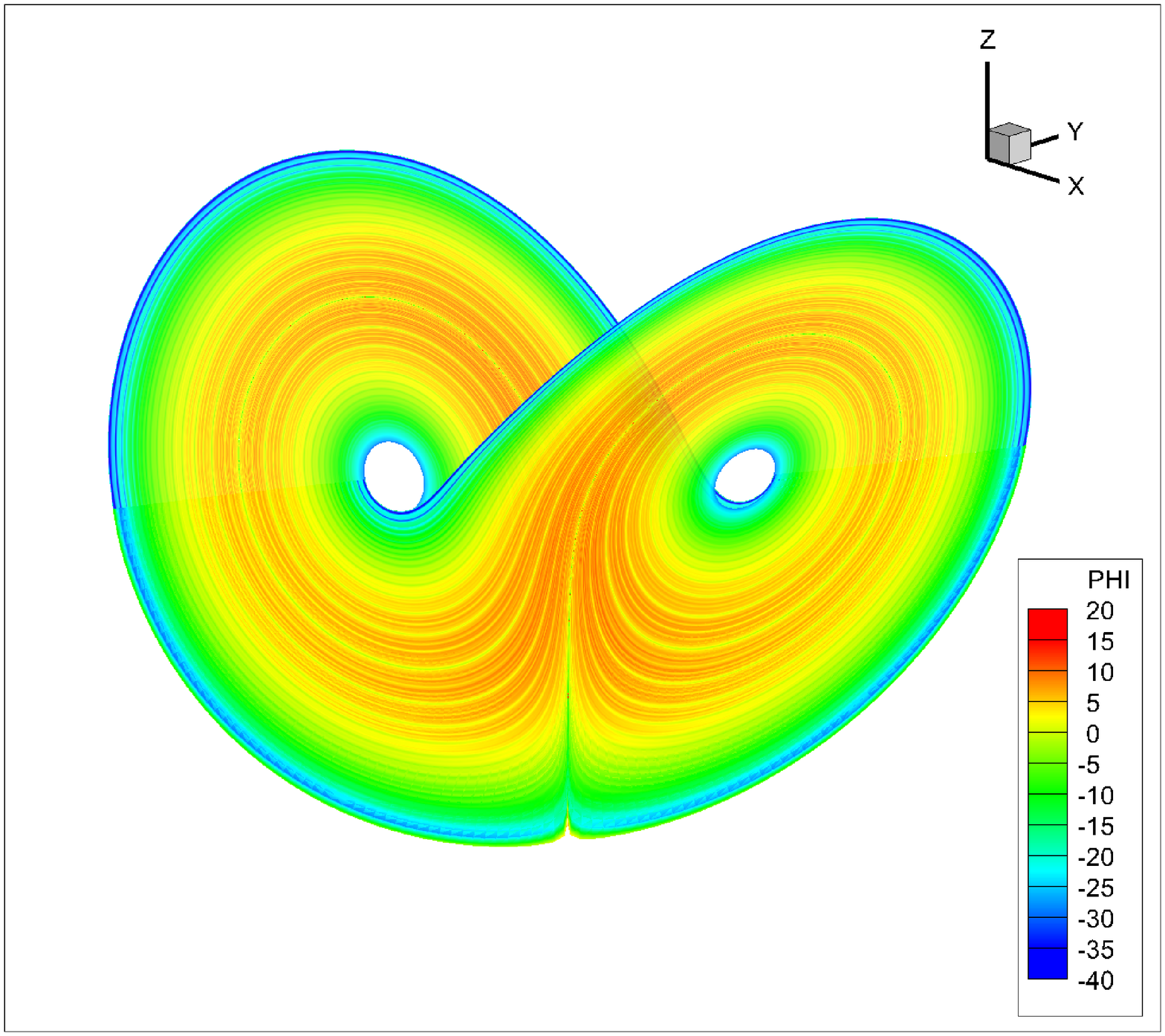}
    \caption{Adjoint distribution on the surface of the Lorenz Attractor for a 
512 by 128 mesh.}
    \label{f:lorenz_phi}
  \end{minipage}
\end{figure}

Figures \ref{f:poincare_phi} and \ref{f:lorenz_phi} show the adjoint density distribution
 on the Poincar\'e section and on the entire attractor surface.  As for the cusp map, 
the adjoint has a fractal structure.  Starting from the Poincar\'e plane, a given 
distribution is duplicated along both branches of the attractor and is propagated backward 
in time towards the origin, where it is squeezed and merged with the distribution from 
the other side of the attractor.  This merged, squeezed distribution then 
propagates back to the Poincar\'e section.  As the sensitivity is for long time 
averages, this process is repeated many times, resulting in the fine fractal 
structures shown in figures \ref{f:poincare_phi} and \ref{f:lorenz_phi}. 

The ability to resolve such fine structures in the adjoint is the main strength of the density adjoint.  Such fine features are in stark contrast to the smooth features of the adjoint distributions computed using other Fokker-Planck approaches, such as those presented by Thuburn \cite{Thuburn:2005:FP}.  The different adjoint solutions arise because many Fokker-Planck approaches, including that of Thuburn, introduce stochasticity to the dynamical system of interest.  This is done to introduce numerical stability, at the cost of reduced accuracy in computing the adjoint and any sensitivities \cite{Thuburn:2005:FP}.  In our method we solve the Fokker-Planck equations on a grid formed using deterministic solutions of the dynamical system.  Most of the numerical dissipation is introduced at the Poincar\'e section, due to the linear interpolation used to form $P_n$.  One interpretation of this is that the invariant measure $\rho_s$ from each streamline can mix at the Poincar\'e section but no where else.  Because of this, the fine 
adjoint structures computed at the Poincare section spread over the entire surface of the attractor manifold, giving us greater insight into the adjoint and sensitivity of chaotic systems than other methods.  

\begin{figure}[htbp]
  \centering

\includegraphics[width=3.2in]{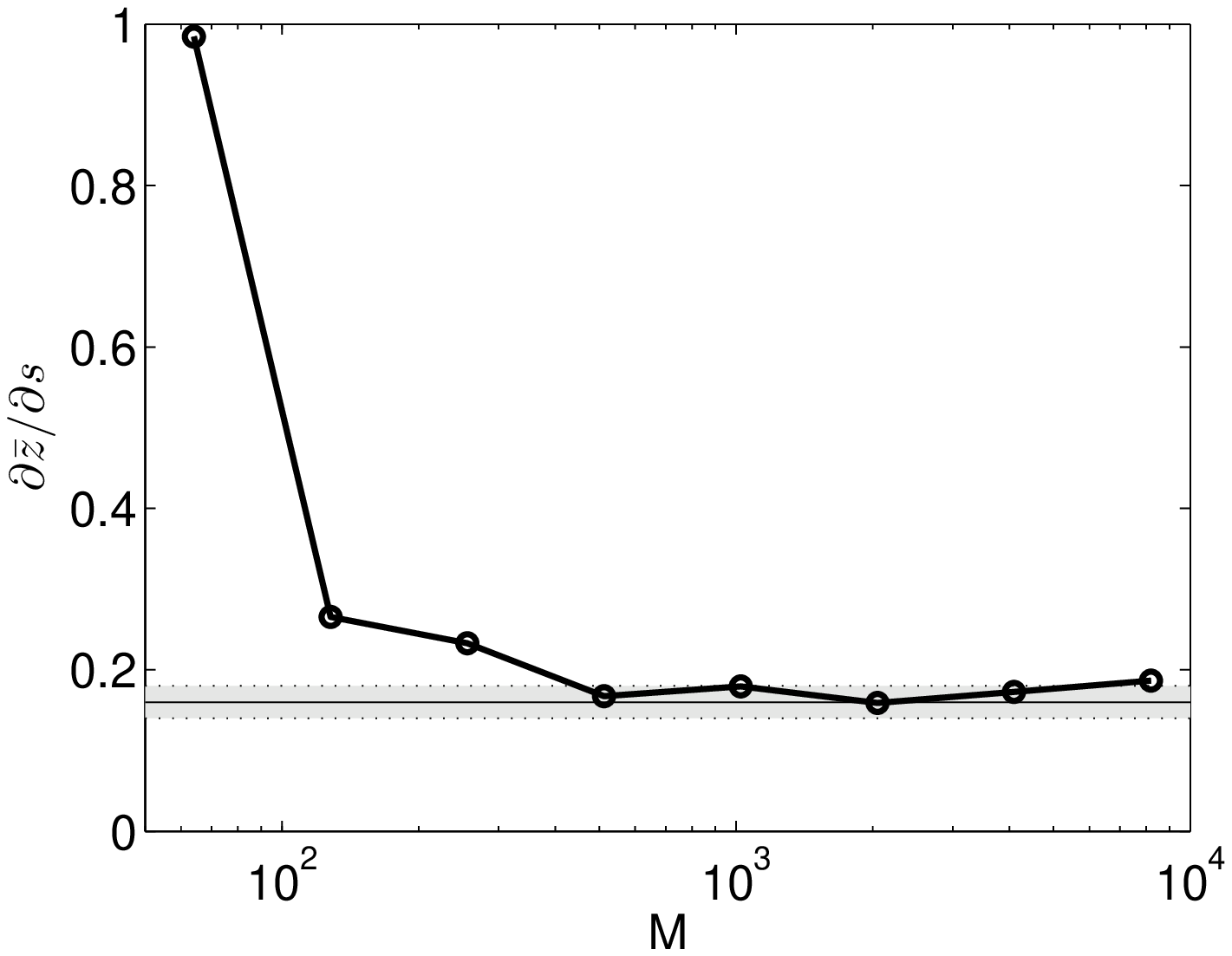}
\includegraphics[width=3.2in]{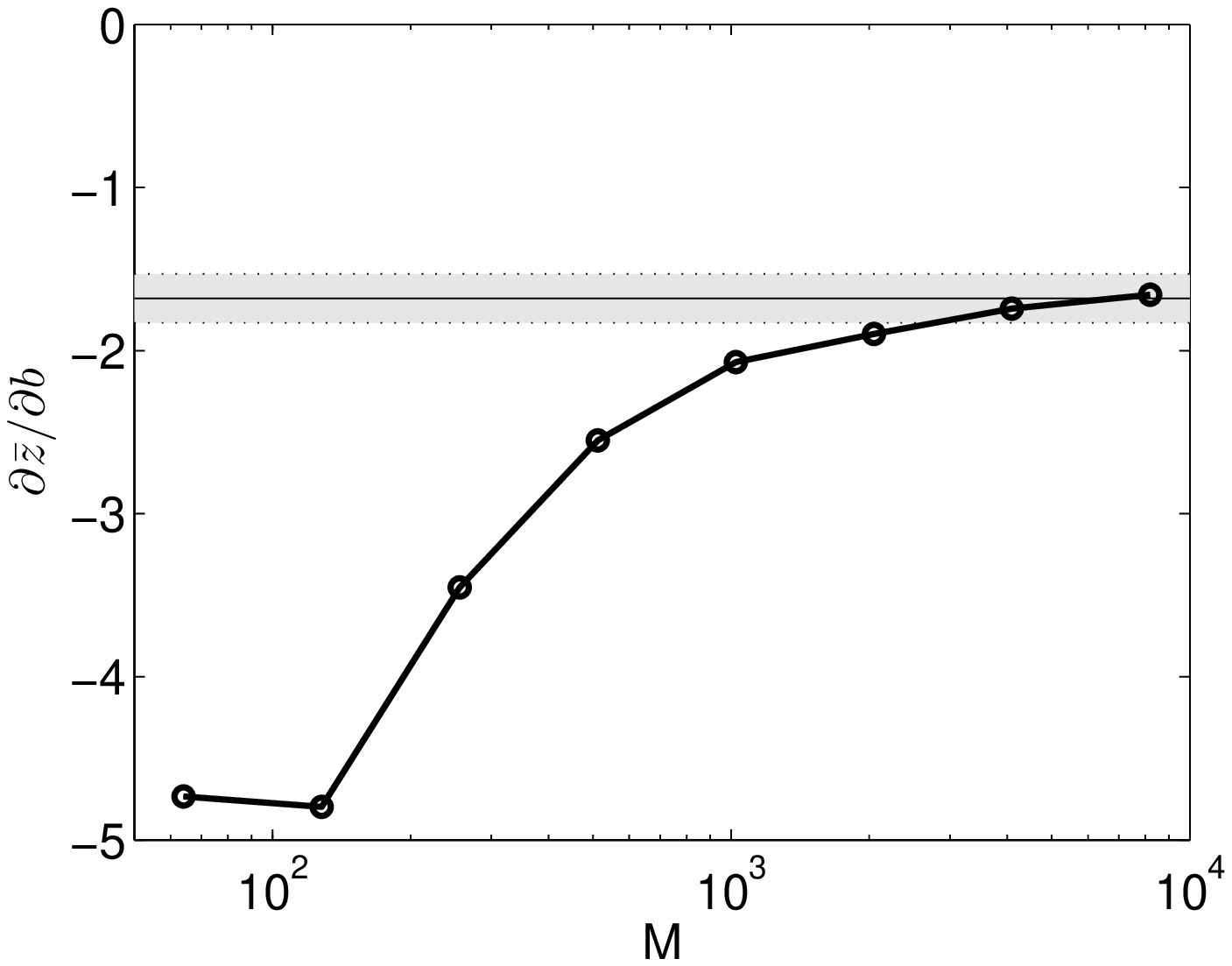}\\
\includegraphics[width=3.2in]{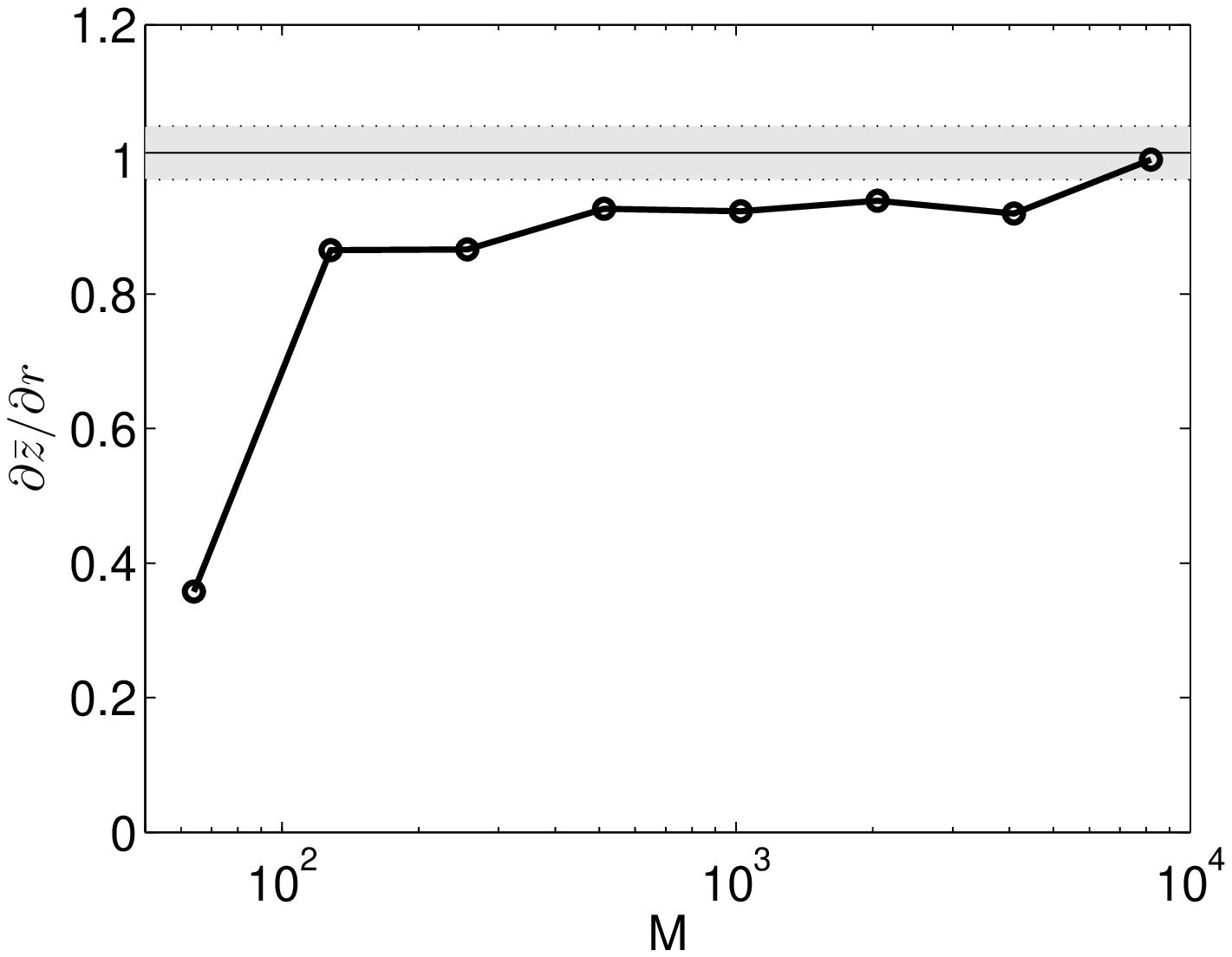}
\includegraphics[width=3.2in]{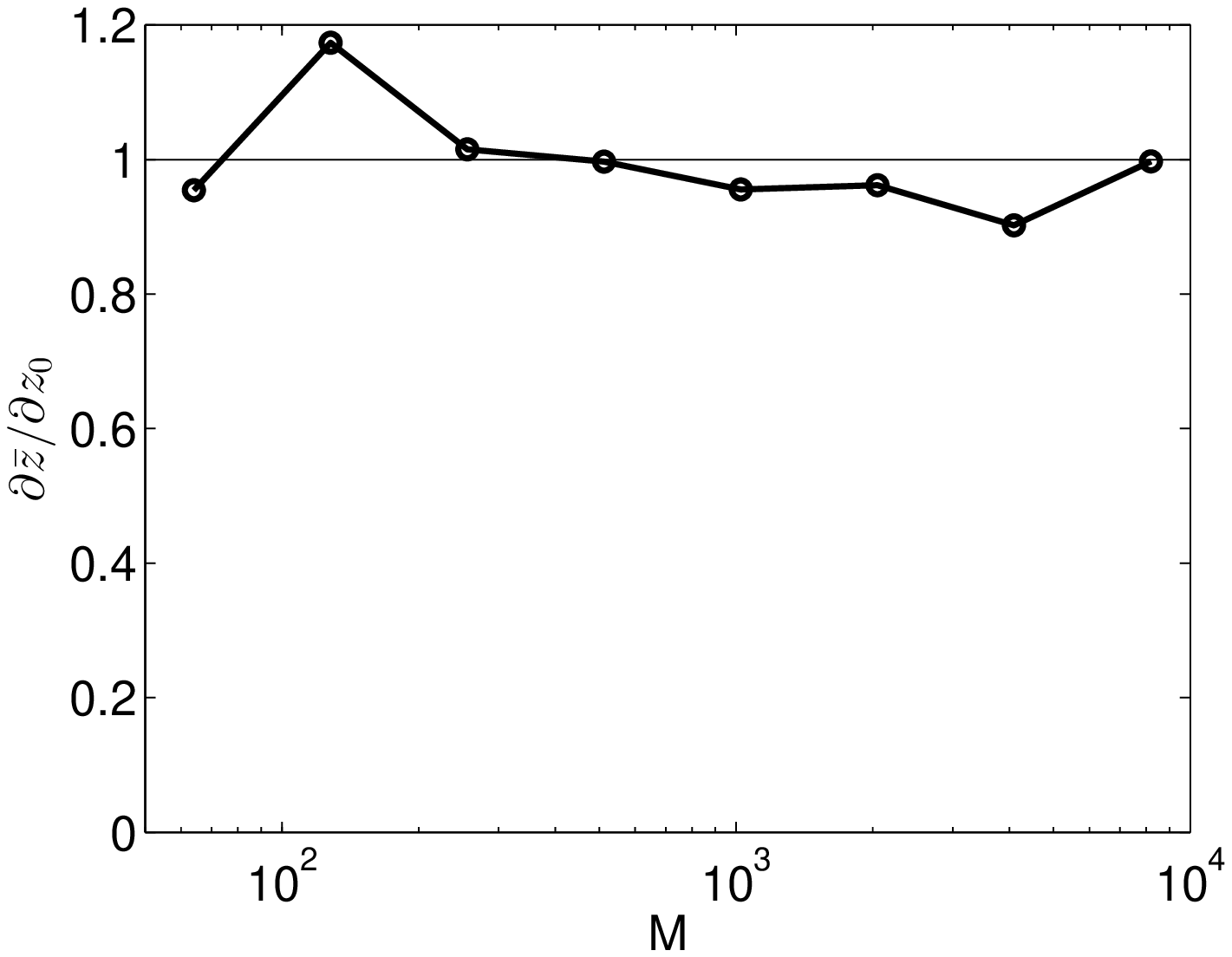}
    \caption{Sensitivity of $\overline{z}$ with respect to the parameters $s$, 
$b$, $r$, and $z_0$ for different numbers of streamlines $M$. The streamlines were distributed so that they were clustered, as in figure \ref{f:lorenz_matrix_NU}.  
$N=128$ is the number of nodes along a streamline. 
The thin black lines correspond to sensitivities computed using a linear regression of ensemble averaged data and the dotted lines are the $3\sigma$ 
confidence bounds from \cite{Wang:2013:LSS1}.  }
    \label{f:lorenz_derivs}

\end{figure}
 
We use the adjoint, in equation \eqref{e:lorenz_grad}, to compute the sensitivity of the average $z$ 
position of the Lorenz attractor with respect to a number of parameters.  From figure \ref{f:lorenz_derivs} it can be seen that the density adjoint 
method predicts the gradients quite well.  Gradients with respect to $s$, $r$ and $b$ parameters are within the $3\sigma$ confidence bounds of the gradients \cite{Wang:2013:LSS1} if the number of streamlines, $M$, is sufficiently large.  Using the highest resolution $M=8192, N =128$ grid the value of $\frac{\partial \overline{z}}{\partial z_0}$ computed was within $0.5\%$ of the correct value of $1.0$.  

\subsubsection*{Sources of Error}

There are a number of sources of error in the density adjoint method, some of which are brought to light by the application of the density adjoint to the Lorenz system.  A key source of error is the discretization of the attractor manifold.  We have approximated the Lorenz attractor as a 2D surface, and while the fractal dimension of the attractor and figure \ref{f:lorenzPoincareMap} indicate that this is a good approximation, it is safe to assume there is still a small amount of error associated with it.  The resolution and structure of the discretization produce a significant portion of the error in computing the gradients.  This can be seen by the sensitivity of gradient computations to the number of streamlines, $M$ in figure \ref{f:lorenz_derivs}.  

 \begin{figure}[htbp]
  \centering
    \begin{minipage}[b]{3.125in}
    \includegraphics[width=3.2in]{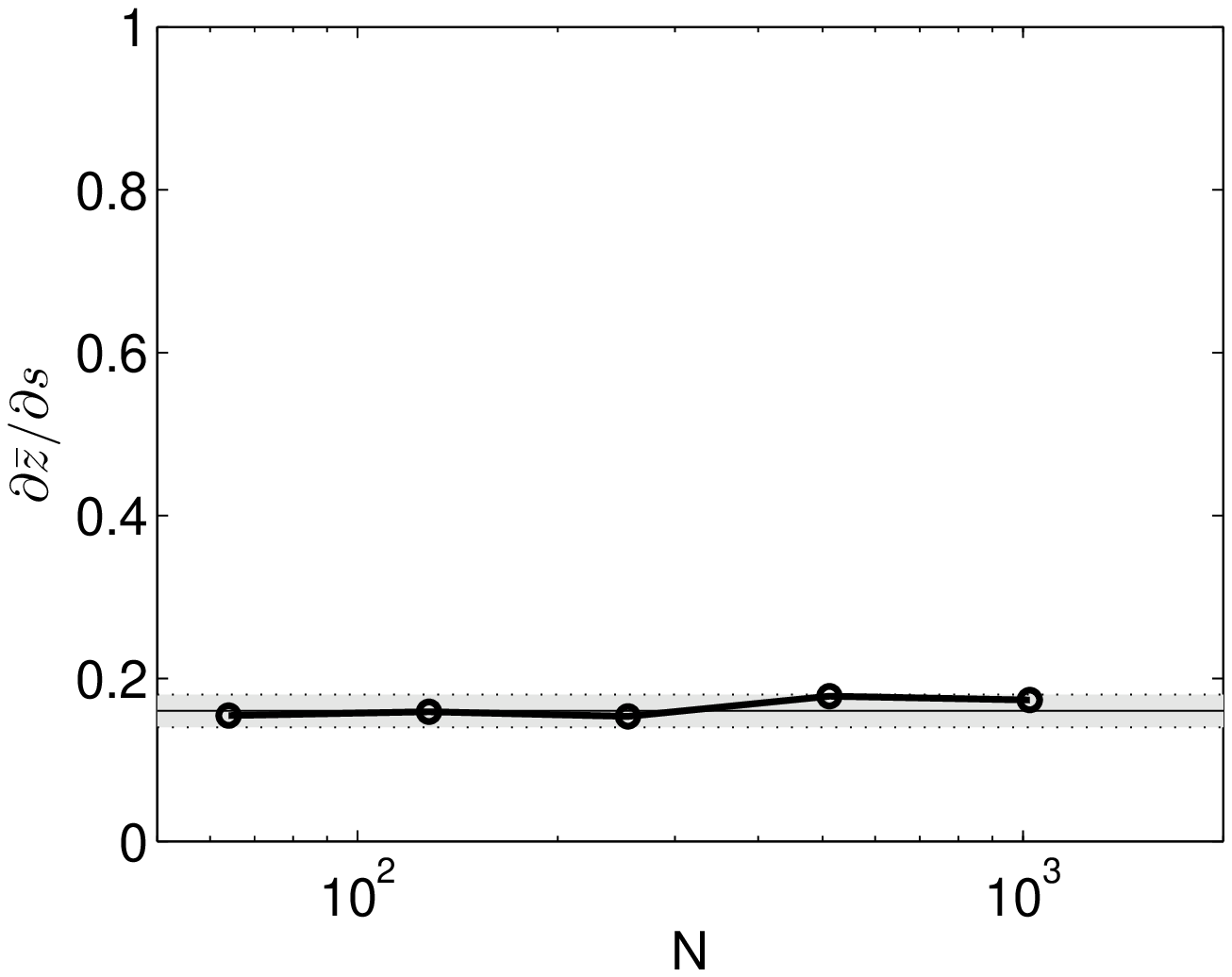}
    \caption{Sensitivity of $\overline{z}$ with respect to the parameter $s$ for different amounts of nodes along each streamline, $N$.  $M=2048$ streamlines were distributed at the Poincare section, as in figure \ref{f:lorenz_matrix}.  The thin black lines correspond to sensitivities computed using a linear regression of ensemble averaged data and the dotted lines are the $3\sigma$ 
confidence bounds from \cite{Wang:2013:LSS1}.  }
    \label{f:ds_N}
  \end{minipage}  
  \hspace{0.125 in}
  \begin{minipage}[b]{3.125in}
    \includegraphics[width=3.2in]{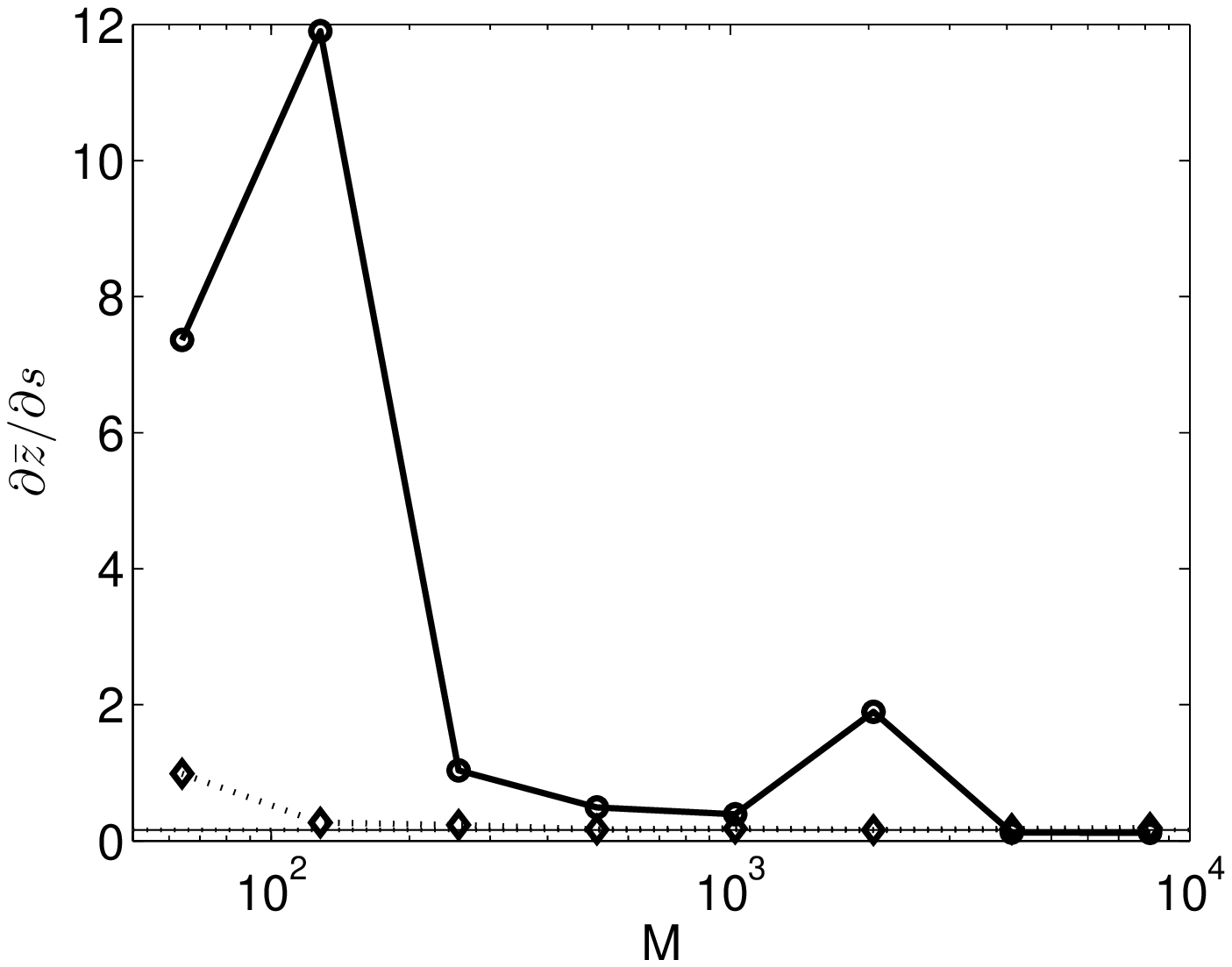}
    \caption{Sensitivity of $\overline{z}$ with respect to the parameter $s$ for different numbers of streamlines $M$.  Convergence for two streamline distributions is shown, uniform $x$-spacing as a solid line and clustered $x$-spacing (from figure \ref{f:lorenz_derivs}) as a dotted line.  
$N=128$ is the number of nodes along a streamline. 
The thin black lines correspond to sensitivities computed using a linear regression of ensemble averaged data and the dotted lines are the $3\sigma$ 
confidence bounds from \cite{Wang:2013:LSS1}.  }
    \label{f:ds_Ugrid}
  \end{minipage}
\end{figure}

On the other hand, it seems that the number of grid points in the streamwise direction, $N$, does not have a great effect on the accuracy of the gradients computed, as indicated by figure \ref{f:ds_N}.  Very similar trends were observed for gradients with respect to $b$, $r$ and $z_0$.  This is because $\phi$ and $\rho_s$ vary more slowly in the $l$ direction than the $s$ direction, as seen in figures \ref{f:lorenz_rho} and \ref{f:lorenz_phi}, respectively.  Overall, it appears that the accuracy of the density adjoint method is more dependent 
on $M$ than on $N$.  

The distribution of the streamlines has a great effect on the convergence of gradient computations.  In figure \ref{f:ds_Ugrid}, we see that the gradient $\partial\overline{z}/\partial s$ converges very slowly when the streamline starting positions are spaced uniformly in $x$ on the Poincar\'e section.  Very similar trends were observed for gradients with respect to $b$, $r$ and $z_0$.  When the streamline starting points are clustered towards $x\approx 13$, much faster convergence is observed. This suggests that the convergence rates of gradients are very sensitive to the discretization of the grid in the spanwise direction.  

A number of numerical techniques and methods used in the density adjoint algorithm (section \ref{ss:alg_ND}).   Proper care with these techniques and methods is needed so that they do not introduce additional numerical errors to the computed gradients.  
Firstly, linear interpolation is used to form the approximate Frobenius-Perron operator $P_n$. One effect of this is that density is not conserved as it flows across the Poincar\'e section, causing the first eigenvalue of $P_n$ to not be exactly $1$. The error due to interpolation can be controlled by using a sufficient number of streamlines to form the mesh.  

Additionally, in step 3 the corresponding eigenvector, $\rho_0$, is found using a power method and then smoothed using a filter.  As long as the threshold frequency of the filter is sufficiently high, gradient accuracy will not be adversely effected.  If the power method is run for a sufficient number of iterations (512 was sufficient for us) it will produce relatively little error.  This assessment is consistent for the gradient results for 1D maps in section \ref{ss:cusp_results}, which are unaffected by attractor displacement and deformation.  

Furthermore, forming the mesh and computing $\rho_s$, $\phi$ and other quantities requires numerical integration, which has some error.  For the Lorenz system, numerical integration in time was conducted with a 4th order Runge-Kutta scheme and a fairly small time step size of $\Delta t = 0.01$.  Using a high order of accuracy and a small time step size will ensure that numerical integration contributes little to errors in gradient computation.  

Some error is also introduced in computing the spanwise direction, $\hat{s}$, gradient of the quantity $\rho_s \frac{\partial \vec{f}}{\partial \xi}$ in equation \eqref{e:lorenz_grad}.  Finite difference approximations are used to compute derivatives  in the streamwise, $\hat{l}$, direction and between the $ith$ streamwise node on adjacent streamlines.  Since the gap between these two nodes is not only in the $\hat{s}$ direction, the $\hat{s}$-direction derivative is approximated by a projection (see appendix \ref{ap:NDgrad} for more details).  However, this error should be reasonably small for a sufficiently smooth $\rho_s \frac{\partial \vec{f}}{\partial \xi}$.  

Finally, in deriving the density adjoint it is assumed that perturbations to parameters of interest (i.e. $s,r,b,$ and $z_0$) can be expressed solely as perturbations to the stationary density $\rho_s$ on the attractor manifold.  However, parameter perturbations can also displace and deform the attractor itself.  consider the Lorenz attractor, which has features that depend on the location of the three fixed points $(\pm\sqrt{b(r-1)},\pm\sqrt{b(r-1)},r+z_0-1)$ and $(0,0,z_0)$.  From inspection, it is clear that the location of these points, which determine the locations of the two holes and the bifurcation, respectively, depend on three of our parameters of interest.  Fortunately the relatively large negative Lyapunov exponent of the Lorenz attractor ensures that attractor deformation and displacement do not contribute much to sensitivities.  However, it might be necessary to take into account the displacement and deformation of the attractor manifold due to parameter perturbations for other chaotic dynamical systems with smaller negative Lyapunov exponents (less dissipative systems). 

\subsubsection*{Computational Costs and Challenges}
Other Fokker-Planck approaches tend to solve the Fokker-Planck equations on a discretization of phase space in the vicinity of the attractor \cite{Thuburn:2005:FP}.  By only solving the Fokker-Planck equation on the attractor manifold, the density adjoint approach has reduced the dimension of the Fokker-Planck equation from that of phase space to that of the attractor manifold.  For the Lorenz system, this means solving a 2D PDE instead of a 3D PDE. However, this does not necessarily result in reduced computational costs.  For example, the coarsest grid used that computed a value of $\partial\overline{z}/\partial s$ within the confidence intervals from \cite{Wang:2013:LSS1} had around 60000 nodes, roughly double the amount used by Thuburn for his Fokker-Planck approach \cite{Thuburn:2005:FP}.  

The size of the grid, especially $M$, is the main driver of the computational cost of the density adjoint method.  $M$ is the number of time integrations required to compute the streamlines needed to form the grid on the attractor manifold in step 2 of the algorithm.  Another $M$ time integrations are required to compute $\rho_s$ and $\phi$ on the attractor surface (step 7). Since these integrations are independent initial value problems, step 2 and step 7 could be carried out in parallel, potentially resulting in a faster solver.  $P_n$ is a sparse $M$ by $M$ matrix, with 4 non-zeros on each row.  Therefore, the cost of the power method in step 3 and the solution of the adjoint matrix system in step 6 both scale with $M^p$, where $p$ is some positive number depending on the solution method.  

Computing gradients, step 9 in the algorithm, has a cost that scales with $MN$, where $N$ is the number of node along each streamline.  However, it is important to note that $N$ did not need to be nearly as large as $M$ to compute accurate gradients, as shown in figure \ref{f:ds_N}.  Additionally it is important to reemphasize that this is an adjoint method.  This means that once steps 1-7 of the algorithm are completed, gradients for any number of parameters can be computed by carrying out steps 8-9 for a given parameter $\xi$.  

However, the previously discussed expenses are incurred after we have obtained a grid approximating the attractor manifold.  Even for the Lorenz system, which has a clear choice for the Poincar\'e plane, a very long time integration of $T=10000$ time units was required to form the Poincar\'e section and conduct the curve fit shown in figure \ref{f:lorenzPoincare2D}.  Also, in order to find the best streamline distribution, the $x$ position on the Poincar\'e section corresponding to the bifurcation at $(0,0,z_0)$ needed to be found, which required a considerable number of numerical experiments.   

Furthermore, there are a number of challenges in applying the density adjoint method that are not touched upon by the Lorenz system example.  Firstly, much of the process behind building the grid approximating the Lorenz attractor depended on being able to visualize the phase space the attractor lies in and the attractor dimension being close to an integer value.  Unfortunately, these two properties do not hold for all chaotic dynamical systems of interest.  Therefore, discretizing attractor manifolds may prove difficult (or impossible) for many chaotic dynamical systems.  Also, even if there were a general method to discretize attractor surfaces, the density adjoint method would be infeasible for higher dimensional attractors.  

Overall, the key limitations of the method are as follows:

\begin{itemize}
\item The use of the Fokker-Planck equation makes the method impractical for high dimensional systems due to high computational costs.  
\item There must be a Poincar\'e section which captures trajectories flowing through the entire attractor.  
\item The intersection of the attractor with the Poincar\'e section must be approximated by a relatively simple curve fit.  
\item Discretizing attractors is relatively straight forward for systems with attractor fractal dimensions that are approximately an integer (i.e. the Lorenz system).  Building a grid on attractor manifolds with non-integer dimensions could prove difficult.
\end{itemize}
 

\section{Conclusion}
\label{s:conclusion}

In conclusion, the density adjoint method computes the sensitivity of 
long time averaged quantities to input parameters for ergodic chaotic dynamical systems if 
a few conditions are met.  Firstly, the system must have a smooth invariant measure.  Secondly, the manifold of the strange attractor must be approximated as an integer-dimensional manifold, which then must be discretized. For 1D chaotic 
maps discretization is trivial.  For continuous chaotic systems, such as the Lorenz 
system, discretization of the attractor manifold is more involved but achievable.  


The density adjoint method computes accurate gradients for the 1D cusp map and the approximately 2D Lorenz attractor. The method also provides insight into adjoint sensitivities of chaotic systems.  The structure of the adjoint solution computed using the density adjoint method is considerably more detailed than that computed used other Fokker-Planck methods.  The adjoint density is observed to be fractal in structure, an illuminating result given that the stretching and folding of density distributions forwards in time becomes compressing and duplicating adjoint density distributions backwards in time. This fractal structure is due to the peak of the cusp map and the bifurcation of the Lorenz attractor and further work is needed to see if fractal adjoints are specific to these systems or if all chaotic dynamical systems have fractal adjoint densities.  The good results obtained using the density adjoint method also shows that accurate gradients can be computed in the presence of some numerical dissipation.  


The density adjoint method could be used to analyze low-dimensional chaotic systems, such as reduced order models of chaotic aero-elastic oscillations of aircraft wings and control surfaces. However the method suffers from the curse of dimensionality like other Fokker-Planck methods.  In addition, the need to discretize the attractor could make it infeasible to extend the density adjoint method to high-dimensional chaotic systems, such as climate models and turbulent aerodynamics simulations. Also, some additional work needs to be done on computing the contribution of the displacement of the attractor manifold to sensitivities to use the method for systems that are less dissipative than the Lorenz system.

%


Despite these limitations, the density adjoint gives more insight into the sensitivity of low dimensional chaotic systems.  The more finely detailed adjoint sensitivity distributions obtained from our method enables better analysis and control of chaotic dynamical systems.


\bibliographystyle{elsarticle-num}
\section*{Bibliography}
\bibliography{main}

\appendix
\include{app1D}

\include{appND}

\end{document}

%% file: app1D.tex
\section{Probability density adjoint for 1D maps}


\subsection{Deriving the continuous density adjoint equation}
\label{ap:1Dcont}

\noindent To derive the adjoint equation, we define a function $v(x) = 1$ and an inner product:

\[
 \langle a, b \rangle = \int_0^1 a(x) b(x) dx
\]

\noindent Consider a small perturbation to $P$.  From the definition of $P$, $(P\rho_s)(x) = \rho_s(x)$, and conservation of 
probability mass:

\begin{equation}
 \delta(P\rho_s) = \delta P\rho_s + P \delta\rho_s = \delta\rho_s, \quad \langle v, \delta \rho_s \rangle = 0
\label{e:density_perturb}
\end{equation}

\noindent Define $\phi$ as the adjoint variable.  Using integration by parts:

\begin{equation}
  0 = \langle \phi, \delta P\rho_s + P \delta\rho_s - \delta\rho_s \rangle = 
\langle P^*\phi - \phi, \delta \rho_s \rangle + \langle \phi, \delta P \rho_s \rangle
\label{e:density_adj1}
\end{equation}

\noindent Combining equations (\ref{e:density_perturb}) and (\ref{e:density_adj1}) with equation (\ref{e:FP_fwd}) : 

\begin{align*}
 \delta \overline{J} &= \langle v, J \rangle - \langle P^*\phi - \phi, \delta \rho_s \rangle 
+ \langle \phi, \delta P \rho_s \rangle - \eta \langle v, \delta \rho_s \rangle = 0 \\
&= \langle J - \eta v - P^*\phi + \phi, \delta \rho_s \rangle + \langle \phi, \delta P \rho_s \rangle
\end{align*}

\noindent For $\phi$ and $\eta$ such that:

\[
 \langle J - \eta v - P^*\phi + \phi, \delta \rho_s \rangle = 0
\]

\noindent Gradients can be computed as follows:  

\begin{equation}
 \delta J = \langle \phi, \delta P \rho_s \rangle
 \label{e:grad2}
\end{equation}

\noindent We derive an expression to compute $\delta P \rho_s$ in section \ref{s:density1d} and \ref{ap:1Dgrad}.  

\noindent To find $\eta$:

\[
 J + \delta J = \langle J, \rho_s + \delta \rho_s \rangle
\]

\noindent For equation (\ref{e:grad2}) to be consistent with this:

\[
 0 = \langle \rho_s, J - \eta v - P^*\phi + \phi \rangle = \langle \rho_s, J \rangle 
- \eta \langle \rho_s, v \rangle - \langle \phi, P\rho_s - \rho_s \rangle
\]

\noindent The second and third inner products on the right hand side are by definition 
1 and 0 respectively, therefore:

\[
 \eta = \langle \rho_s, J \rangle = \overline{J}
\]

\noindent Therefore, the adjoint equation is:

\[
 P^*\phi - \phi = \overline{J} - J
\]

\subsection{Derivation of the gradient equation}
\label{ap:1Dgrad}

%
%
%
%
%
%
%
%
%

\begin{figure}[htb!] \centering
\includegraphics[width=3.2in]{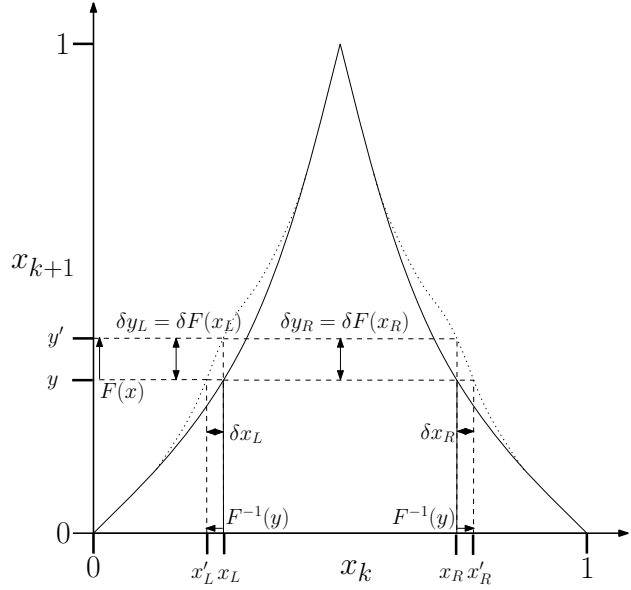}
\caption{The effect of a perturbation on the mapping function.  }
\label{f:delta_f}
\end{figure}

From figure \ref{f:delta_f}, a local functional perturbation left of the peak 
moves $F^{-1}(y)$ by $\delta x_L$ to the left.  Also, note that $\delta F/\delta 
x_L$ is a first order approximation to the local slope, therefore $\delta F = 
-F'(F^{-1}(y))\delta x_L$ for small $\delta F$.  A similar argument can be made to show that the same 
equation applies to the right of of the peak.  Therefore equation \ref{e:cdf_perturb2} 
can be rewritten as:

\[
 \int_0^y \delta \rho_0 \ ds = \frac{\rho_s(x_L)}{F'(x_L)}\delta F(x_L) - \frac{\rho_s(x_R)}{F'(x_R)} \delta 
F(x_R)
\]

%

\subsection{Deriving the discrete density adjoint equation}
\label{ap:1Ddisc}

\noindent Recall that the first eigenvalue of the discrete operator $P_n$ is not exactly
one.  Denoting $\lambda$ as the first eigenvalue, which converges to one as $n \to \infty$:

\[
 P_n \hr_s - \lambda \hr_s = 0
\]

\noindent Now consider the linearization:

\[
 \delta P_n \hr_s - \delta \lambda \hr_s + P_n \delta \hr_s - \lambda 
\delta \hr_s = 0, \quad -\hv^T \delta \hr_s = 0
\]

\noindent Combine this with the discrete version of equation \ref{e:FP_fwd}:

\[
 \delta \overline{J} = \frac{1}{n} [ J^T \delta \hr_s + \hp^T (  
- \delta \lambda \hr_s + P_n \delta \hr_s - \lambda \delta \hr_s) 
+ \hp^T \delta P_n \hr_s - \eta v^T \delta \hr_s ]
\]

\noindent Where $\eta$ is the adjoint variable for the eigenvalue perturbation 
$\delta \lambda$. Rearrange to isolate $\delta \hr_s$ and $\delta \lambda$:

\[
 \delta \overline{J} = \frac{1}{n} [\hp^T \delta P_n \hr_s + 
(J^T - \hp^T P_n - \hp^T \lambda - \eta \hv^T )\delta \hr_s + 
(-\hp^T \hr_s)\delta \lambda]
\]

\noindent The matrix expression in section \ref{s:density1d} is obtained by eliminating any dependence of
$\delta \overline{J}$ on $\delta \hr_s$ and $\delta \lambda$.  Also, from the continuous adjoint
equation, it can be seen that $\eta = \overline{J}$:

\begin{align*}
 P_n^T \hp  - \lambda \hp - \hv \overline{J} &= J \\
-\hr^T_s \hp &= 0 
\end{align*}

%% file: appND.tex
\section{Probability density adjoint for continuous chaos}

\subsection{Divergence operator on the attractor manifold}
\label{ap:div}

\begin{figure}[htb!]
  \centering 
    \includegraphics[width=2.5in]{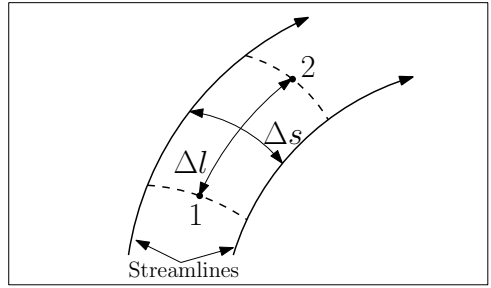}
    \caption{Schematic of two streamlines on the attractor manifold.  }
    \label{f:div_scheme}
\end{figure}

\noindent The divergence operator can be derived by considering the weak form of the operator:

\begin{equation}
 \int_{\Omega} \vec{\nabla}_s \cdot \vec{f} \ dA = \int_{\partial \Omega} f \cdot \hat{n} \ dS
 \label{e:weak_div}
\end{equation}

Where $\Omega$ is the area enclosed between the two streamlines and the boundaries 1 and 2 shown in figure \ref{f:div_scheme}, $\partial \Omega$ is the boundary of this area, and $\hat{n}$ is the unit normal vector pointing outwards from the boundary.  

Since $\vec{f}$ is tangent to the streamlines by definition, only the boundaries at 1 and 2 contribute to the integral on the right hand side of equation \eqref{e:weak_div}.  The boundaries at 1 and 2 can be chosen so that $\hat{n} = \hat{l}$.  If this is the case, equation \eqref{e:weak_div} can be rewritten as follows for an infinitesimal area $\Omega$:

\begin{equation}
 \vec{\nabla}_s \cdot \vec{f} \Delta l \Delta s = (\vec{f}(l_2,s)\cdot \hat{l})(\Delta s + \delta s) - (\vec{f}(l_1,s)\cdot \hat{l})\Delta s
 \label{e:div_inf}
\end{equation}

\noindent Where $l_1 = l$, $l_2 = l + \Delta l$, $s_1 = s_2 = s$, and the lengths of the boundaries 1 and 2 are $\Delta s$ and $\Delta s + \delta s$, respectively.  

\noindent We can find an expression for $\delta s$ by considering the linearized governing equations:

\begin{equation}
 \frac{d\vec{\delta x}}{dt} = \textbf{J} \vec{\delta x}
\end{equation}

\noindent Where $\textbf{J}$ is the Jacobian of $\vec{f}(\vec{x})$.  The solution, $\vec{\delta x}$, is a first order approximation for the separation between two adjacent streamlines.  For our infinitesimal area $\Omega$:

\[
 \frac{\vec{\delta x}_2 - \vec{\delta x}_1}{\Delta t} \approx \textbf{J}\vec{\delta x}_1 
\]

\noindent Since $\vec{\delta x_1} \approx (\Delta s) \hat{s}$, it follows that $ \vec{\delta x}_2  \approx \textbf{J}\hat{s} \Delta s \Delta t + (\Delta s) \hat{s}$, and since the $\hat{s}$ component of $\vec{\delta x}_2$ is approximately $\Delta s + \delta s$:

\[
\delta s \approx (\textbf{J}\hat{s})\cdot \hat{s} \Delta s \Delta t = (\hat{s}^T \textbf{J} \hat{s}) \Delta t \Delta s
\]

\noindent Using the above expression, along with the Taylor expansion $\vec{f}(l+\Delta l,s) = \vec{f}(l,s) + \textbf{J} \cdot \hat{l} \Delta l$, equation \eqref{e:div_inf} becomes:

\begin{equation}
 \vec{\nabla}_s \cdot \vec{f} \Delta l \Delta s = (\hat{l}^T \textbf{J} \hat{l})\Delta l \Delta s + (\hat{s}^T \textbf{J} \hat{s}) (\vec{f}(l,s)\cdot \hat{l}) \Delta t \Delta s
 \label{e:div_inf2}
\end{equation}

Substituting the identity $ (\vec{f}(l,s)\cdot \hat{l}) \Delta t = |\vec{f}(l,s)| \Delta t = \Delta l$ and dividing equation \eqref{e:div_inf2} by $ \Delta l \Delta s$, we obtain an expression for the divergence operator in terms of the Jacobian, $\textbf{J}$ and the unit vectors $\hat{l}$ and $\hat{s}$:

\[
 \vec{\nabla_s}\cdot \vec{f} = \hat{l}^T \textbf{J} \hat{l} + \hat{s}^T \textbf{J} \hat{s}
\]

\subsection{Deriving the continuous adjoint equation}
\label{ap:NDcont}

\noindent First, linearize equation (\ref{e:dens_cons}):

\begin{equation}
\vec{\nabla}_s \cdot (\delta \rho \vec{f} + \rho \delta \vec{f}) = 0
\label{e:dens_perturb}
\end{equation}

As (\ref{e:dens_perturb}) is zero, it can be multiplied by some scalar variable 
$\phi$ and added to equation (\ref{e:delta_J}).  By conservation of probability mass, a 
perturbation to $\rho_s$ does not change the total probability:

\[
\iint \delta \rho \ dl ds = 0
\]

\noindent Therefore, $\delta \rho$ can also be added to equation (\ref{e:delta_J}):

\begin{equation}
\delta\overline{J} = \iint \phi \vec{\nabla}_s \cdot (\delta \rho \vec{f} + \rho_s \delta 
\vec{f}) + J(\vec{x})\delta \rho + c\delta \rho \ dl ds
\label{e:delta_J_adj}
\end{equation}

\noindent Where $c$ is some constant. Conducting integration by parts:

\begin{align*}
\delta\overline{J} &= \iint \phi \nabla_s \cdot (\delta \rho \vf )+ 
J(\vec{x})\delta \rho  + c\delta \rho + \phi \nabla_s \cdot (\rho_s \delta \vf) 
\ dl ds \\
&= \iint -\delta \rho \vf \cdot \nabla_s \phi + J(\vec{x})\delta \rho + c\delta 
\rho + \phi \nabla_s \cdot (\rho_s \delta \vf)) \ dl ds \\
&= \iint (-\frac{\partial \phi}{\partial t} + J(\vec{x}) + c)\delta \rho + \phi 
\nabla_s \cdot (\rho_s \delta \vf)) \ dl ds 
\end{align*}


\noindent In order to eliminate the dependence of $\delta\overline{J}$ on $\delta \rho$, the adjoint density equation is:

\begin{equation}
\frac{\partial \phi}{\partial t} = J(\vec{x}) + c
\label{e:lorenz_adj1}
\end{equation}

\noindent Multiplying both sides of (\ref{e:lorenz_adj1}) by $\rho_s$ and 
integrating over the attractor surface shows that $c = -\overline{J}$:

\begin{align*}
\iint \rho_s f(\vec{x}) \nabla_s \phi \ dl ds &= \iint \rho_s J(\vec{x}) + \rho_s 
c \ dl ds\\
\iint \phi \nabla_s \cdot (\rho_s f(\vec{x})) \ dl ds &= \iint \rho_s J(\vec{x}) 
+ \rho_s c \ dl ds\\
0 &= \iint \rho_s J(\vec{x})\,dl ds +  c\iint \rho_s \ dl ds\\
c &= -\overline{J} \\
\end{align*}

Therefore:

\[
\frac{\partial \phi}{\partial t} = J(\vec{x}) -\overline{J} 
\]

If the above equation is satisfied, equation (\ref{e:delta_J_adj}) reduces to

\[
 \delta \overline{J} = \iint \phi \, \vec{\nabla}_s \cdot (\rho_s \delta \vec{f}) \ dl ds
\]

\subsection{Deriving the discrete Adjoint Equation}
\label{ap:NDdisc}

\noindent Consider the eigenvalue equation 
for the Poincar\'e stationary density:

\[
P_n \hr_0 = \lambda \hr_0
\]
 
\noindent This can be modified using $D$, a diagonal matrix with $|\vec{f}(\vec{x}_0) \times \hat{s}_0| ds_0$ 
for the ith streamline along the main diagonal:

\[
D P_n D^{-1} D \hr_0 = \lambda D \hr_0
\]

\noindent Defining $D \hr_0 = \hq$ and $D P_n D^{-1} = A$:

\[
A \hq = \lambda \hq
\]

The adjoint is derived using $A$ and $\hq$ because perturbations to $\hq$ 
correspond to density perturbations on the attractor surface.  A 
perturbation $\delta \hq$ can be written as follows:

\begin{align}
\delta (\lambda \hq) &= \delta (A \hq) \nonumber\\
\lambda \delta \hq + \hq \delta \lambda &=  \delta A \hq + A \delta \hq \nonumber\\
(\lambda I - A) \delta \hq - \delta A \hq + \hq \delta \lambda &= 0
\label{e:delta_q_zero}
\end{align}

\noindent From equation (\ref{e:avg_frm_rho0}), $\delta \overline{J}$ is related to a perturbation 
to $\hq$ as follows

\begin{equation}
\delta \overline{J} = \hJ^T \delta \hq
\label{e:disc_grad_q}
\end{equation}

\noindent where we define $\hJ_i = \int_0^T J(t) dt$ for streamline $i$.  

\noindent Also, it can be shown that:

\begin{equation}
\iint \delta \rho_s dl ds \quad \Rightarrow \quad \hv^T D^{-1} \delta \hq = 0
\label{e:q_cons}
\end{equation}

\noindent Adding equation (\ref{e:disc_grad_q}) to the product of equation (\ref{e:delta_q_zero}) and 
the discrete density adjoint $\hp_0$ as well as equation (\ref{e:q_cons}) and the adjoint eigenvalue 
$\eta$ yields:

\begin{align*}
\delta \overline{J} &= \hJ^T \delta \hq + \hp_0^T ( (\lambda I - A) \delta \hq + \eta(\hv^T D^{-1} \delta \hq)
-\delta A \hq + \hq \delta \lambda  ) \\
\delta \overline{J} &= (\hp_0^T (\lambda I - A) + \eta \hv^T D^{-1} + \hJ^T)  \delta \hq + \hp_0^T \hq \delta \lambda  - \hp_0^T 
\delta A \hq  
\end{align*}

\noindent To eliminate the dependence of $\delta \overline{J}$ on $\delta \hq$ and 
$\delta \lambda$:

\[
(A-I)^T \hp_0 + \eta D^{-1} \hv = \hJ , \qquad  \hq^T \hp_0
\]

\noindent Therefore:

\[
(DP_nD^{-1}-I)^T \hp_0 + \eta D^{-1} \hv = \hJ , \qquad  \hr_s^T D \hp_0
\]
As in section \ref{s:density1d}, it can be shown that 
$\eta = -\overline{J}$, therefore:

\[
\left[ \begin{array}{cc} (D^{-1}P^T D-\lambda I) & -D^{-1} \hv \\ \hr_s^T D & 
0  \end{array} \right] \left[ \begin{array}{c} \hp_0 \\ \overline{J}  
\end{array} \right] = \left[ \begin{array}{c} \hJ \\ 0 \end{array} 
\right]
\]


\subsection{Computing Attractor Surface Areas}
\label{ap:attractor_area}

\noindent Recall:

\[
\rho(\vec{x}) |\vec{f}(\vec{x})| ds = \rho_0 |\vec{f}_0 \times \hat{s}_0| ds_0
\]

Where $ds$ is the streamline width and the subscript 0 indicates values at the 
start of the streamline.  $ds_0$ is set as the average distance from a 
streamline to its neighboring streamlines along the Poincar\'e section.  

Noting that $|\vec{f}(\vec{x})| = \frac{\partial l}{\partial t}$:

\[
|\vec{f}(\vec{x})| ds = \frac{\partial l}{\partial t} ds = \frac{\partial A}{\partial 
t}
\]

Therefore:
\begin{equation}
\frac{\partial A}{\partial t} = \frac{\rho_0 |\vec{f}_0 \times \hat{s}_0| 
ds_0}{\rho(\vec{x})}
\label{e:area}
\end{equation}

Integrating this equation along a streamline yields the total area of that 
streamline.  For a given node $k$, $dA_k$ is found by taking the difference of $A$ 
at the midpoint between nodes $k-1$ and $k$ and the midpoint between nodes $k$ and 
$k+1$.  The area of the first node is computed as $A$ at the first midpoint.  The 
area of the last node is the difference between $A$ for the entire streamline and 
the last midpoint.  

\subsection{Computing gradients on the attractor surface}
\label{ap:NDgrad}

The partial derivative in the $l$ direction is found using a central difference 
when possible.  At the beginning of a given streamline, a forward difference is 
used and a backward difference is used at the end of a given streamline. As 
$\rho_s \frac{\partial \vec{f}}{\partial \xi}$ is a three dimensional vector, three 
derivatives are obtained, corresponding to the x, y and z components.  

To find the partial derivative in the $s$ direction, first consider the forward 
difference between the $j$th node on streamline $i$ and the $j$th node on streamline 
$i+1$.  This difference can be used to approximate the derivative in the $s'$ 
direction, which is not equal to the $s$ direction.  To find the $s$ direction, 
the difference of some vector $\vec{X}$ in the $s'$ direction, $\Delta_{s'} \vec{X}$  
can be decomposed as follows:

\[
\Delta_{s'} \vec{X} = \alpha \frac{\Delta_{l} \vec{X}}{\Delta l} + \beta \frac{\Delta_{s} 
\vec{X}}{\Delta s}
\]  

Where $\vec{X} = \rho_s \frac{\partial \vec{f}}{\partial \xi}$, $\alpha = \vec{s'} \cdot 
\hat{l}$ and $\beta = \vec{s'} \cdot \hat{s}$.  It is important to note that 
$\vec{s'}$ is not a unit vector like $\hat{l}$ and $\hat{s}$.  This expression 
can be rearranged to yield an expression for $\frac{\partial \vec{X}}{\partial s}$:

\[
\frac{\partial \vec{X}}{\partial s} \approx \frac{\Delta_{s} \vec{X}}{\Delta s} = 
\frac{1}{\beta} \Delta_{s'} \vec{X} - \frac{\alpha}{\beta} \frac{\Delta_{l} \vec{X}}{\Delta 
l}
\]

This same equation can be solved for the backwards difference and the average 
of the forward and backward differences can be taken to find the central 
difference.  

\noindent Finally, to find the surface gradient:

\[
\nabla_s \cdot \left(\rho_s \frac{\partial f}{\partial \xi}\right) = 
\frac{\partial \vec{X}}{\partial l} \cdot \hat{l} + \frac{\partial \vec{X}}{\partial s} 
\cdot \hat{s}
\]

\noindent Therefore:

\[
\frac{\partial \overline{J}}{\partial \xi} \approx \sum_{k=0}^{N} \phi_k \, 
\left[\frac{\partial \vec{X}}{\partial l} \cdot \hat{l} + \frac{\partial \vec{X}}{\partial 
s} \cdot \hat{s}\right]_k dA_k
\]